\documentclass{amsart}
\usepackage{amssymb}
\usepackage{url}
\usepackage{amscd}
\usepackage{adjustbox}
\theoremstyle{cupthm}

\newtheorem{theorem}{Theorem}[section]
\newtheorem{proposition}[theorem]{Proposition}
\newtheorem{corollary}[theorem]{Corollary}
\newtheorem{definition}[theorem]{Definition}
\newtheorem{question}[theorem]{Question}
\newtheorem{remark}[theorem]{Remark}

\newtheorem{example}[theorem]{Example}

\usepackage{graphicx}
\usepackage {tikz}

\def\cp{$\frak{g}=(g_1,\dots, g_n)$ and $\mathcal
E=\{E_1,\dots, E_m\}$ be a finite string of elements of $G$ and a
finite partition for $G$, respectively. }

\def\c{Con(\frak g,\mathcal E)}

\usetikzlibrary {positioning}
\usepackage {xcolor}

\definecolor {processblue}{cmyk}{0.96,0,0,0}
\usetikzlibrary{arrows, decorations.markings}

\tikzstyle{vecArrow} = [thick, decoration={markings,mark=at position
   1 with {\arrow[semithick]{open triangle 60}}},
   double distance=1.4pt, shorten >= 5.5pt,
   preaction = {decorate},
   postaction = {draw,line width=1.4pt, white,shorten >= 4.5pt}]
\tikzstyle{innerWhite} = [semithick, white,line width=1.4pt, shorten >= 4.5pt]
\usetikzlibrary{intersections}
\usepackage{tkz-euclide}
\usetkzobj{all}

\begin{document}

\markboth{Akram Yousofzadeh}
{A constructive way to compute the Tarski number of a group}

%
%

\title{A CONSTRUCTIVE WAY TO COMPUTE THE TARSKI  NUMBER OF A GROUP}

\author{\footnotesize AKRAM YOUSOFZADEH
}


\email{ayousofzade@yahoo.com}

\maketitle


\begin{abstract}
The Tarski number of a group $G$ is the minimal number of the pieces of paradoxical decompositions of that
 group. Using configurations along with a matrix combinatorial property we construct paradoxical decompositions.
 We also compute an upper bound for the Tarski number of a given non-amenable group by counting the number of paths
  in a diagram associated to the group.
\end{abstract}



\section{Introduction} \label{sect1}
 Rosenblatt and Willis introduced a concept for
groups  to show that for an infinite discrete amenable group or a
non-discrete amenable group $G$ a net of positive, normalized
functions in $L^1(G)$ can be constructed such that this net
converges weak* to invariance but does not converge strongly to
invariance \cite{r-w}. This concept which is called {\it configuration}
is also used to classify some
group theoretical properties (see for example \cite{arw, ary}).

Configurations are strongly linked to the amenability of groups and by Tarski's
alternative, a discrete group is
non-amenable if and only if it admits a paradoxical decomposition.
Therefore it is valuable to construct the paradoxical
decomposition for such a group, using configurations. This problem
which was originally asked by Willis is answered partially in
\cite{y-r}. In that paper the paradoxical decomposition was
constructed under a paradoxical condition.

In the present  paper we define a general matrix combinatorial
property under which the paradoxical decomposition is completely
constructed. We also find a new upper bound for the  Tarski number of
a given non-amenable group.

\subsection*{Notations}

 The following notations are used throughout this paper\\

\begin{itemize}
\item $\mathbb{N}$, $\mathbb{Z}$ and
$\mathbb{R}$ are the sets of natural, integer and real numbers,
respectively,

\item $\bigsqcup$ is the disjoint union of sets,

\item $gA=\{ga;\ \ a\in A\},$  for a group $G$, $A\subseteq G$
and $g\in G,$

\item $\mathcal P(X)$ is the power set of the set $X$,

\item $|X|$ is the the cardinal number of the set $X$,

\item A $(0,1)$-matrix is a matrix with entries in $\{0,1\}$.

\end{itemize}

If $v$ is a vector of real numbers and $n\in \mathbb{Z}$, by $v\geq n$ we
mean that $v$ has entries equal to or greater than n.

\section{Preliminaries} \label{ns}
\subsection{Matrix theory}
We use Gordan's theorem to find nontrivial nonnegative solutions to  a
homogenous system of linear equations with real coefficients. This theorem  has also
applications in linear programming \cite{linear programming}.

\begin{theorem}(Gordan 1873)\label{gordan}
Either a linear homogenous system of equations $AX=0$ possesses a
nontrivial solution in nonnegative variables or there exists an
equation, formed by taking some linear combination of equations,
that all positive coefficients. That is, either there exists an
$x$ such that $$Ax=0,\ \ \ 0\neq x\geq 0$$ or there exists a
vector $m$ such that $m^tA>0$ (has positive entries).
\end{theorem}
\begin{remark}\rm In Theorem \ref{gordan}, if $A$ is a matrix with rational entries,
then the entries of $m$ can be chosen in $\mathbb{Z}$.
\end{remark}
 The main theorem of this
paper is proved under a certain condition. To clarify this condition, we need some
definitions.

Let $\pi:\{1,\dots,n\}\rightarrow \{1,\dots,n\}$ be a permutation
for the set $\{1,\dots,n\}$. Then $$P_{\pi}=\left(
\begin{array}{c}
e_{\pi(1)} \\
e_{\pi(2)} \\
\vdots \\
e_{\pi(n)}
\end{array} \right)$$ is
called the permutation matrix associated to $\pi$, where $e_{i}$
denotes the row vector of length $n$ with 1 in the i-th position
and 0 otherwise. When the permutation matrix $P_{\pi}$ is
multiplied with a matrix $M$ from left, $P_{\pi}M$ will permute
the rows of $M$ by $\pi$.

If $P=\left( \begin{array}{cccc}
P_1 \\
P_2 \\
\vdots \\
P_n
\end{array} \right)$ is a permutation matrix, by $P^{+}$ we mean the matrix with
shifted rows, i.e. $$P^{+}=P_{\rho}P=\left( \begin{array}{cccc}
P_2 \\
P_3 \\
\vdots \\
P_n\\
P_1
\end{array} \right),$$
in which $\rho$ is the cyclic permutation $(1\ 2\ \dots \ n)$.
 Throughout we use the notation $$T=\left(
\begin{array}{ccccc}
1& 0& 0& \dots & 0 \\
1& 1& 0& \dots& 0 \\
\vdots \\
1& 1& 1& \dots& 1
\end{array} \right)=\sum_{1\leq j\leq i\leq n}E_{ij},$$
where $E_{ij}$ is the matrix with 1 in ${ij}$ position and 0
otherwise. When the matrix $T$ is multiplied with a matrix $M$
from left, j-th row of $TM$ will be the sum of j first rows of
$M$.

\begin{definition}\rm
Let $\ell\in \mathbb{N}$ and $\left( \begin{array}{cccc}
A_1 \\
A_2 \\
\vdots \\
A_n
\end{array} \right)
$ and $  \left( \begin{array}{cccc}
B_1 \\
B_2 \\
\vdots \\
B_n
\end{array} \right)
$ be two $(0,1)$-matrices with rows $A_i, B_i$. Let also the
vector $\sum_{i=1}^n (B_i-A_i)$ has strictly positive entries. If
there exists a permutation matrix $P_{\pi}$ such that the matrix
\begin{equation}\label{normality}
TP_{\pi}(B-A)-P_{\pi}^{+}A
\end{equation}
has integer entries equal or greater than $-1$, we say that the
homogenous system of equations $(B-A)X=0$ is {\it{normal.}}
\end{definition}

It is apparent that if $\sum_{i=1}^n (B_i-A_i)$ has strictly positive entries,
then the system has no non-zero non-negative solution. Conversely,
if $(B-A)X=0$ is a system of equations with no non-zero
non-negative solution, then by Gordan's theorem  there exists a
vector $m=(m_1,\dots, m_n)$ such that $m^t(B-A)$ has strictly
positive entries. If we permit $B-A$ to have repeated rows and
we insert the opposite of a row (exchanging the corresponding rows of
$A$ and $B$) if necessary, then $m_i$ can be chosen in $\{0,1\}.$
Omit $B_i-A_i,$ where $m_i=0$ and denote the modified matrix by
$B-A$ again. So, we can assume that $\sum_{i=1}^n
(B_i-A_i)$ has strictly positive entries.

It is to be noted that there exist examples of both
normal and non-normal systems of equations with no nonzero non-negative solution.
We will explain  the relation between normality and paradoxical decompositions in section 3.

\subsection{Non-amenable discrete groups}
Let $G$ be a discrete group. Then $G$ is called amenable if it
admits a finitely additive probability measure $\mu$ on the
$\sigma-$algebra $\mathcal P(G)$ such that
$$\mu(gA)=\mu(A),\ \ \ \ \ \ \ (A\subseteq G,\ g\in G).$$

\begin{definition} \rm\cite{wagon} \label{norm}
Let $G$ be a group acting on a set $X$ and suppose $E\subseteq X$.
$E$ is $G$-paradoxical (or, paradoxical with respect to $G$) if
for some positive integers $m, n$ there are pairwise disjoint
subsets $A_1,\dots,A_n,B_1,\dots,B_m$ of $E$ and
$g_1\dots,g_n,h_1\dots,h_m\in G$ such that
\begin{equation}\label{star}
 E =\bigsqcup_{i=1}^ng_iA_i=\bigsqcup_{j=1}^mh_jB_i.\end{equation}
\end{definition}
A group $G$ is called paradoxical if it is $G$-paradoxical, where
$G$ acts on itself by left multiplication. Clearly if $G$ is a
paradoxical group satisfying the above definition, then it cannot
be amenable. Indeed if $\mu$ is a $G-$invariant probability
measure, then \begin{align*}1&=\mu(G)=\sum\mu(A_i)+\sum
\mu(B_j)\\&=\sum\mu(g_iA_i)+\sum\mu(h_jB_j)=\mu(G)+\mu(G)=2.\end{align*}
In fact there is the following remarkable alternative due to
Alfred Tarski.
\begin{theorem}
Let $G$ be a discrete group. Exactly one of the following
happens

1)   $G$ is paradoxical,

2)  $G$ is amenable.

\end{theorem}

 The number $\tau=n+m$ for $n$ and $m$ in (\ref{star}) is called the Tarski number of that
paradoxical decomposition; the minimum of all such numbers over
all the possible paradoxical decompositions of $G$, is called the
Tarski number of $G$ and denoted by $\tau(G)$. In the case that
there is no paradoxical decomposition, we set $\tau(G)=\infty$.
The  Tarski number of a group is of real interest and has been
estimated for some classes of groups. But it is not so
easy to compute in many cases.
For more details see \cite{sapir} and
\cite{silber}.

 We draw the reader's attention to the next proposition for different types of paradoxical decompositions.

\begin{proposition}\cite[Proposition 1.2]{y-r}\label{belg}
Let $G$ be a group. Then the following statements are equivalent

1)  There exist a partition $\{A_1,\dots,A_n,B_1,\dots,B_m\}$ of
$G$ and $g_1,\dots,g_n$ and $h_1,\dots,h_m$ in $G$ such that
$\{g_iA_i\}_{i=1}^n$ and $\{h_jB_j\}_{j=1}^m$ form partitions of
$G$.

2)  There exist pairwise disjoint subsets $A_1,\dots,A_n,B_1,\dots,B_m$ of $G$ and
elements $g_1,\dots,g_n$ and $h_1,\dots,h_m$ in $G$ such that
$\{g_iA_i\}_{i=1}^n$ and $\{h_jB_j\}_{j=1}^m$ form partitions of
$G$.

3)  There exist pairwise disjoint subsets
$A_1,\dots,A_n,B_1,\dots,B_m$ of $G$ and elements $g_1,\dots,g_n$
and $h_1,\dots,h_m$ in $G$ such that $G =\bigcup_{i=1}^n
g_iA_i=\bigcup_{j=1}^m h_jB_j$ (not necessarily pairwise
disjoint).
\end{proposition}

Because of the above equivalence, $G$ admits a
paradoxical decomposition, if any condition of Proposition
\ref{belg} holds. The decomposition in condition 1 is
called a \emph{complete} paradoxical decomposition.
By \cite[Theorem A.1]{sapir ershov} if the conditions 2 or 3 are satisfied
for a group $G$, then $\tau(G)\leq m+n$.

\subsection{Configuration of groups}
Let $G$ be a discrete group. The configurations of $G$ are defined
in terms of finite generating sets and finite partitions of $G$.
If $\frak g=(g_1,\dots,g_n)$ is a string of elements of $G$ and
$\mathcal E = \{E_1,\dots, E_m\}$ is a partition of $G$, a
configuration corresponding to $(\frak g,\mathcal E)$  is an
$(n+1)-$tuple $C=(c_0,\dots,c_n)$, where $1\leq c_i\leq m$ for
each $i$, such that there is $x$ in $G$ with $x \in E_{c_0}$ and
$g_ix\in E_{c_i}$ for each $1\leq i\leq n$. The set of all
configurations corresponding to the pair $(\frak g,\mathcal E)$
will be denoted by $\c$. It is shown that groups with the same set
of configurations have some common properties. For example they
obey the same semigroup laws and have the same Tarski numbers (see
\cite{arw} and \cite{ytr}).

In the case that $\frak g=\{g_1,\dots,g_n\}$ is a generating set
for $G$,  the configuration $C =(c_0,\dots,c_n)$ may be described
as a labelled tree which is a subgraph of the Cayley graph of the
finitely generated group $G$ and configuration set $\c$ is a set
of rooted trees having height 1. In last section of the paper we
assign a new graph to $G$ that depends on the pair $(\frak
g,\mathcal E)$.

 If  $(\frak g,\mathcal E)$ is as above and for each $C\in \c$
$$x_0(C)=E_{c_0}\cap(\cap_{j=1}^n g_j^{-1}E_{c_j})\ \ \ \ \text{and}\ \ \ \ x_j(C)=g_jx_0(C),$$
  then it is seen that for any $0\leq j\leq n,$
$\{x_j(C);\ \ C\in\c \}$ is a partition for $G$.
Let $C\in \c$ and $f\in \ell^1(G).$ Define $f_C=\sum_{x\in x_0(C)}f(x).$ Then we have (see \cite{r-w})
$$\langle f-_{g_j}f,\chi_{E_i} \rangle=0,\ \ \ \ \ \ (1\leq j\leq n,\ 1\leq i\leq m)$$ if and only if
$$\sum_{x_0(C)\subseteq E_i}f_C=\sum_{x_j(C)\subseteq E_i}f_C\ \ \ \ \ \ (1\leq j\leq n,\ 1\leq i\leq m).$$
For each pair $(\frak g,\mathcal E)$
for $G$, the system of equations
$$\sum_{x_j(C)\subseteq E_i}f_C=\sum_{x_k(C)\subseteq E_i}f_C,\ \ \ \  \ (1\leq i\leq m,\ 0\leq j,k\leq n) $$
with variables $f_C,\ \ C\in \c$  is called the system of
configuration equations corresponding to $(\frak g,\mathcal E)$
and is denoted by $Eq(\frak g,\mathcal E)$.
By a normalized solution to this system, we mean a solution $(f_C)_C$ such that for each $C$, $f_C\geq 0$ and $\sum_{C}f_C=1$.
A group $G$ is amenable if and only if there is a net $(f_{\alpha})_{\alpha}$ of positive functions in $L^1(G)$ with $\int f_{\alpha}=1$ which is
weak* convergence to invariance, that is for any $F\in L^{\infty}(G)$ and any $g\in G$
$$\lim_{\alpha}\langle f_{\alpha}-_{g}f_{\alpha},F \rangle=0$$
(see \cite{paterson}). Using this fact, Rosenblatt and
Willis proved the following theorem
\begin{theorem}\cite[Proposition 2.4]{r-w}\label{roswill}
There is a normalized solution of every possible instance of the
configuration equations if and only if $G$ is amenable.
\end{theorem}
\noindent We will apply this theorem to find paradoxical decompositions.

\section{Main Theorem}
Throughout this section by $G$, $\frak{g}=(g_1,\dots,
g_n)$ and $\mathcal E=\{E_1,\dots, E_m\}$ we mean a group, a finite string of
elements of $G$ and a finite partition for $G$, respectively.
If $D$ is a subset of $Con(\frak g,\mathcal E)$, we use the
following notation
$$\tilde{D}:=\bigsqcup_{C\in D}x_0(C).$$
In particular $\tilde{\c}=G$. Clearly  $D_1,D_2\subseteq \c$ are disjoint if and only if $\tilde D_1$ and $\tilde D_2$ are. The
configuration equation $\sum_{x_j(C)\subseteq
E_i}f_C=\sum_{x_k(C)\subseteq E_i}f_C$ is written in the form
$\frak aX=\frak bX,$ where $\c=\{C_1,\dots, C_{\ell}\}$, $$ X=
\left(\begin{array}{c}f_{C_1}\\f_{C_2}\\ \vdots\\f_{C_{\ell}}
\end{array} \right),
$$ $\frak a$ is the
coefficient vector of the left hand side and $\frak b$ is the
coefficient vector of the right hand side of the equation.

\begin{theorem}\label{main}
 If a subsystem of $Eq(\frak g,\mathcal E)$ is normal, then $G$
is non-amenable and a paradoxical decomposition of $G$ can be
written in terms of $\frak g$ and $\mathcal E$.
\end{theorem}
\begin{proof}
Let  $V=\left( \begin{array}{cccc}
V_1 \\
V_2 \\
\vdots \\
V_p
\end{array} \right)$ and  $W=\left( \begin{array}{cccc}
W_1 \\
W_2 \\
\vdots \\
W_p
\end{array} \right)$ be $(0,1)$-matrices in $ M_{p\times\ell}(\mathbb{R})$
such that $(W-V)X=0$ is the normal subsystem of $Eq(\frak
g,\mathcal E)$ satisfying (\ref{normality})for a permutation matrix $P_{\pi}$. Write
$$V_t=(V_t(C))_{C\in\c}\ \ \ \text{and}\ \ \ \ W_t=(W_t(C))_{C\in\c}, \ \ \ \ \ (1\leq t\leq p).$$
Then by the normality of $(W-V)X=0$, for each $C\in\c$
\begin{equation}\label{positivity}
 \alpha_C:=\sum_{t=1}^pW_t(C)-\sum_{t=1}^pV_t(C)>0.
\end{equation}
We consider only the case that $P_{\pi}=I,$ the identity matrix. Otherwise defining  $V^*_t=V_{\pi(t)}$ and $W^*_t=W_{\pi(t)}$ we see that
$$TP_{\pi}(W-V)-P_{\pi}^+V= T(W^*-V^*)-{V^*}^+$$
and the proof is carried out in a similar way.
Therefore throughout the proof we assume
\begin{align*}
W_1(C)-(V_1(C)+V_2(C))\geq& -1,\\
(W_1(C)+W_2(C))-(V_1(C)+V_2(C)+V_3(C))\geq& -1,\\
&\vdots\\
(\sum_{t=1}^{p-1}W_t(C))-\sum_{t=1}^{p}V_t(C))\geq& -1.
\end{align*}
Suppose that  $({i_1},{i_2},\dots,{i_p})$,
$(j_{i_1},j_{i_2},\dots,j_{i_p})$ and
$(k_{i_1},k_{i_2},\dots,k_{i_p})$ are strings such that ${i_s}\in
\{1,\dots,m\}$ and $j_{i_s}, k_{i_s}\in \{0,1,\dots,n\}$ and the
modified system is
$$\sum_{x_{j_{i_t}}(C)\subseteq
E_{i_t}}f_C=\sum_{x_{k_{i_t}}(C)\subseteq E_{i_t}}f_C,\ \ \ \ \ \ \ \
(1\leq t\leq p).
$$
Note that the strings are used instead of subsets since the
repetition is not excluded for the equations.
 For convenience  we use the following notations for $1\leq t\leq p$
$$A_t=\{C;\ \ x_{k_{i_t}}(C)\subseteq E_{i_t}\}\  \text{and}\  B_t=\{C;\ \
x_{j_{i_t}}(C)\subseteq E_{i_t}\}, \ \ \ \ \ (1\leq t\leq p).$$ In other words, the system
can be written as
\begin{equation}\label{ss}
\sum_{C\in A_t}f_C=\sum_{C\in B_t}f_C,\ \ \ \ \ \ \ \
(1\leq t\leq p).
\end{equation}
It is not difficult to see that for $1\leq t\leq p$
\begin{equation}\label{key}
\dot{g}_t \bigsqcup_{C^*\in A_t} x_0(C^*)=\bigsqcup_{C\in B_t} x_0(C),
\end{equation}
where $\dot{g}_t= g_{k_{i_t}}^{-1}g_{j_{i_t}}$. For $C\in B_t$ and $C^*\in A_t$ put $$ E_t(C^*,C)=x_0(C^*)\cap \dot{g}_t^{-1}x_0(C).$$ Then we have by (\ref{key})
\begin{equation}\label{salib}
x_0(C)=\bigsqcup_{C^*\in A_t}\dot{g}_tE_t(C^*,C),\ \ \ \ \ \ \ \ \ \ \ (C\in B_t).
\end{equation}\\

\noindent {\bf Case 1. {The sets $A_t$ are pairwise disjoint}}.
According to (\ref{positivity}) and (\ref{ss}), for each $C\in \frak A:=\bigsqcup_{t=1}^p A_t$ there exist at least two numbers $1\leq t_1^C, t_2^C\leq p$ such that $C\in B_{t_1^C}$ and $C\in B_{t_2^C}$. Also for each $C\in\frak B:= \c\setminus \frak A$ there exists at least one $1\leq t^C\leq p$  such that $C\in B_{t^C}$. This fact together with (\ref{salib}) implies that
$$x_0(C)=\bigsqcup_{C^*\in A_{t_1^C}}\dot{g}_{t_1^C}E_{t_1^C}(C^*,C)=\bigsqcup_{C^*\in A_{t_2^C}}\dot{g}_{t_2^C}E_{t_2^C}(C^*,C),\ \ \ \ \ (C\in \frak A)$$ and $$x_0(C)=\bigsqcup_{C^*\in A_{t^C}}\dot{g}_{t^C}E_{t^C}(C^*,C),\ \ \ \ \ (C\in \frak B).$$ Consequently the fact that
  $$G = \left[\bigsqcup_{C\in \frak A}x_0(C)\right]\bigsqcup\left[\bigsqcup_{C\in \frak B}x_0(C)\right]$$ necessitates the following two equations
  \begin{align}\label{11}
 G=  \left[\bigsqcup_{C\in \frak A}\bigsqcup_{C^*\in A_{t_1^C}} \dot{g}_{t_1^C}E_{t_1^C}(C^*,C)\right]\bigsqcup\left[\bigsqcup_{C\in \frak B}x_0(C)\right]
\end{align}
and
\begin{align}\label{12}
  G= \left[\bigsqcup_{C\in \frak A}\bigsqcup_{C^*\in A_{t_2^C}}\dot{g}_{t_2^C} E_{t_2^C}(C^*,C)\right]\bigsqcup\left[\bigsqcup_{C\in \frak B}\bigsqcup_{C^*\in A_{t^C} }\dot{g}_{t^C} E_{t^C}(C^*,C)\right].
\end{align}
We emphasize that  all the sets in the following families are pairwise disjoint
\begin{enumerate}
  \item $x_0(C),$ with $C\in \frak B$,
  \item $E_{t^C}(C^*,C),$ with $C\in \frak B, C^*\in A_{t^C}$,
  \item $E_{t_1^C}(C^*,C),$ with $C\in \frak A, C^*\in A_{t_1^C}$,
  \item $E_{t_2^C}(C^*,C),$ with $C\in \frak A, C^*\in A_{t_2^C}$.
\end{enumerate}
We know that $\frak A\cap \frak B=\emptyset$ and $A_i\cap A_j=\emptyset$ for $i\neq j.$ On the other hand for $C^*\in A_t,$ $E_{t}(C^*,C)\subseteq x_0(C^*)\subseteq \tilde A_t\subseteq \tilde{\frak A}$.
 Therefore the sets of types (2), (3) and (4) are all disjoint from the sets of type (1). Since $t_1^C\neq t_2^C$, $E_{t_1^C}(C^*,C)\subseteq A_{t_1^C}$ and $E_{t_2^C}(C^*,C)\subseteq A_{t_2^C}$ and $A_{t_1^C}\cap A_{t_2^C}=\emptyset$, we have $E_{t_1^C}(C^*,C)\cap E_{t_2^C}(C^*,C)=\emptyset.$ Now let $i\in\{1,2\}$ and $t_i^{C'}=t^C$, for $C\in \frak B$ and $C'\in \frak A$.
Let $C^*\in A_{t^C}$ and ${C'}^*\in A_{t_i^{C'}}$. Applying  $\frak A\cap \frak B=\emptyset$ once again, we see $g^{-1}_{t^C}x_0(C)\cap g^{-1}_{t_i^{C'}}x_0(C')=\emptyset.$ Thus  $E_{t_i^{C'}}(C'^*,C')\cap E_{t^C}(C^*,C)=\emptyset.$ Therefore all the sets of types (3) and (4) are disjoint from the sets of type (2). It remains to show that for distinct configurations $C,D\in \frak A,$ $i,j\in\{1,2\},$ $C^*\in A_{t_i^C}$ and $D^*\in A_{t_j^D}$ one has  $E_{t_i^C}(C^*,C)\cap E_{t_j^D}(D^*,D)=\emptyset.$ But it is clear since if $t_i^C\neq t_j^D$, then
$$E_{t_i^C}(C^*,C)\cap E_{t_j^D}(D^*,D)\subseteq x_0(C^*)\cap x_0(D^*)\subseteq \tilde A_{t_i^C}\cap \tilde A_{t_j^D}=\emptyset$$
and if  $t_i^C=t_j^D$, then
$$E_{t_i^C}(C^*,C)\cap E_{t_j^D}(D^*,D)\subseteq \dot{g}^{-1}_{t_i^C} (x_0(C)\cap x_0(D))=\dot{g}^{-1}_{t_i^C}(\emptyset)=\emptyset.$$
Consequently the equations (\ref{11}) and (\ref{12}) form a paradoxical decomposition of $G$.\\

\noindent{\bf Case 2. {The sets $A_t$ are not pairwise disjoint}}.
The procedure in case 1 does not work here because the sets in a paradoxical decomposition must be disjoint. In the sequel we shall replace the sets $A_1,\dots,A_p$ with new sets $P_1,\dots,P_p,$ which are disjoint and then a paradoxical decomposition with respect to a partition of $\bigsqcup_{i=1}^p\tilde{P}_i$ will be given.

For $C\in \c$ if $t_0$ is the least number in $\{1,\dots,p\}$ such that $V_{t_0}(C)\neq 0$, define $V'_{t_0}(C)=1$ and for $t\neq t_0$ put $V'_{t}(C)=0$. Then \begin{equation}\label{vt}
\sum_{t=1}^pV'_t(C)=\begin{cases}
                        1, & \mbox{if } \sum_{t=1}^{p}V_t(C)>0 \\
                        0, & \mbox{otherwise}.
                      \end{cases}\end{equation}   We need a counter to compute the number of appearance of a configuration in $V_1,\dots,V_p$ and  $W_1,\dots,W_p$. To achieve this purpose define   \begin{equation}\label{'} m^0_C=\begin{cases}
                       (\sum_{t=1}^{p}V_t(C))-1, & \mbox{if } \sum_{t=1}^{p}V_t(C)>0 \\
                       0, & \mbox{otherwise}
                     \end{cases} \end{equation}
and for $1\leq t\leq p$
 \begin{equation}\label{36'}
    m^{t}_C=\begin{cases}
                    m^{t-1}_C-1, & \mbox{if } W_{t}(C)=1,\ m^{t-1}_C>0 \\
                    m^{t-1}_C, & \mbox{otherwise}.
                  \end{cases}
                  \end{equation}
Consequently by (\ref{'}) we have for $C\in \c$
\begin{equation}\label{mv}
  \sum_{t=1}^pV_t(C)-\sum_{t=1}^pV'_t(C)=m^0_C.
\end{equation}
Consider the vectors $(W'_t(C))_{\c},$ $1\leq t\leq p$ with the following definition
$$W'_{t}(C)=\begin{cases}
   W_{t}(C)-1, & \mbox{if } W_{t}(C)=1,\ m^{t-1}_C>0 \\
  W_{t}(C), & \mbox{otherwise}.
 \end{cases}$$
It is seen by (\ref{36'}) that for every $C\in\c$
\begin{eqnarray*}
  \sum_{t=1}^pW_t(C)-\sum_{t=1}^pW'_t(C)&=& \sum_{t=1}^p[W_t(C)-W'_t(C)] \\
   &=&\sum_{t=1}^p[m^{t-1}_C - m^{t}_C]
    \\
   &=& m^{0}_C-m^{p}_C \\
   &=& m^0_C.
\end{eqnarray*} Thus by (\ref{mv}) and (\ref{positivity})
\begin{equation}\label{subtract}
  \sum_{t=1}^pW'_t(C)-\sum_{t=1}^pV'_t(C)=\sum_{t=1}^pW_t(C)-\sum_{t=1}^pV_t(C)=\alpha_C.
\end{equation}
Put $$P_t:=\{C:\ \ \ V'_t(C)\neq 0\}\ \ \text{and}\ \ \ Q_t:=\{C:\ \ \ W'_t(C)\neq 0\}.$$
Then by (\ref{vt}) it is evident that $P_1,\dots,P_p$ are pairwise disjoint. It is also clear  for $1\leq t\leq p$ that $P_t\subseteq A_t$, $Q_t\subseteq B_t$ and
$A_t\setminus P_t=\bigcup_{j=1}^{t-1}(A_t\cap A_j).$
Let $1\leq t<j\leq p$ and define
 \begin{equation}\label{def d} D_t^{j}=\left[({B}_t\setminus \bigsqcup_{s=t+1}^{j-1} D_t^s)\cap({A}_{j}\setminus {P}_{j})\right]\setminus\bigsqcup_{k=1}^{t-1} D_k^j ,\end{equation} using the conventions $\bigcup_{l=l_1}^{l_2} E_l=\emptyset$, when $0\leq l_2<l_1$.
Obviously
$$D_t^j\cap D_t^{j'}=\emptyset\ \ \ \ \text{and}\ \ \ \ D_t^j\cap D_{t'}^{j}=\emptyset,\ \ \ \ \ \ (t\neq t', j\neq j').$$
Next for $1\leq t\leq p-1$ we prove that \begin{equation}\label{q_t}
                   Q_t=B_t\setminus(\bigsqcup_{j=t+1}^p D_t^j).
                   \end{equation}
 To see this, let $C\in Q_t.$ Then $C\in B_t$ and $W'_t(C)=1,$ which means that $m_C^{t-1}=0$. If $m_C^{0}=0,$ then for each $2\leq j\leq p$, $C\notin A_j\setminus P_j.$
  In particular $C\notin D_t^j,$ for every $t+1\leq j\leq p.$ Assume that $m_C^{t-1}=k>0.$ Then by (\ref{'}) there are exactly $k+1$ numbers $j_0<j_1<\dots<j_{k}$ such that $C\in\cap_{i=0}^k A_{j_i}$ or equivalently
   $$\ \ \ \ \ \ \ \ \ \ \ \ \ \ \ \ \ \ C\in A_{j_i}\setminus P_{j_i},\  \ \ \ \ \ \ \ \ \  (1\leq i\leq k).$$ This along with the normality condition ($\sum_{j=1}^{t-1}W_j-\sum_{j=1}^{t}V_j\geq -1$) implies that there exist natural numbers
    $t_1<t_2<\dots<t_k<t$ such that $C\in B_{t_i},\ (1\leq i\leq k),$ where $t_1$ is the smallest number with this property. Therefore by
    (\ref{def d}) and the sentence after that $C\in \bigcap_{i=1}^k D_{t_i}^{j_i}.$ So,  we have $C\notin \bigsqcup_{j=t+1}^p D_t^j.$ The converse is obtained by a similar argument.
 We continue the proof of the theorem in the following four steps.\\
\\
\noindent {\bf{Step 1.}}
By the normality condition if $1\leq t\leq p$, then $\sum_{j=1}^{t-1}W_j-\sum_{j=1}^{t}V_j\geq -1$. This means  $(A_t\setminus P_t)\subseteq \cup_{j=1}^{t-1}B_j$. We shall prove that for $1\leq t\leq p$
 \begin{equation}\label{b2}{A}_t\setminus {P}_t=\bigsqcup_{s=1}^{t-1}  D_s^t.\end{equation}
Firstly we show that $$ {A}_t\setminus {P}_t\subseteq  \bigcup_{s=1}^{t-1}\left[({B}_s\setminus \cup_{j=s+1}^{t-1}D_s^j)\setminus\bigsqcup_{k=1}^{s-1} D_k^t \right].$$
Let $C\in {A}_t\setminus {P}_t$ and pick the smallest natural number $t_0<t$ such that $C\in B_{t_0}$ (it is possible by the normality condition). If for each $t_0<j<t$ and $1\leq k<t_0$, $C\notin D_{t_0}^j$ and $C\notin D_{k}^j$, we are done. If $C\in \bigsqcup_{k=1}^{t_0-1} D_k^t$, then $C\in \bigsqcup_{k=1}^{t-1} D_k^t$ because $t_0<t.$ Assume that for some $t_0<j<t$, $C\in D_{t_0}^j.$ Then by (\ref{def d}) $C\in {A}_{j}\setminus {P}_{j}.$ By the definition of $P_j$, this means that there exists at least one natural number $s<j$ such that $C\in A_j\cap A_s$ but on the other side $C\in A_t$. Consequently $C\in A_j\cap A_s\cap A_t.$ Since by the normality condition $\sum_{j=1}^{t-1}W_j-\sum_{j=1}^{t}V_j\geq -1$, there are at least two sets of the form $B_i$ containing $C$. Accordingly there exists the smallest positive integer $t'_0<t$ such that $C\in B_{t'_0}.$ This process is obviously finite and therefore we achieve the purpose.  Hence
 \begin{align*}
   {A}_t\setminus {P}_t=({A}_t\setminus {P}_t)\cap \bigcup_{s=1}^{t-1}\left[({B}_s\setminus \cup_{j=s+1}^{t-1}D_s^j)\setminus\bigsqcup_{k=1}^{s-1} D_k^t \right]=\bigsqcup_{s=1}^{t-1}  D_s^t.
 \end{align*}
\noindent{\bf{Step 2.}} We show by induction that for all $j\geq 1$ there exists a family   $\{\mathcal A_{\sigma_j},\ \ \ \sigma_j\in \Lambda_j\}$ with cardinality $2^{j-1}-1$ of disjoint subsets (possibly the empty sets) of $\bigsqcup_{t=1}^{j-1}\tilde{P}_t$ and the members  $\{g_{\sigma_j},\ \ \ \sigma_j\in \Lambda_j\}$ of $G$ such that
 $$\tilde{A}_j\setminus \tilde{P}_j=\bigsqcup_{\sigma_j\in \Lambda_j}g_{\sigma_j}{\mathcal A}_{\sigma_j}.$$
Since  $\tilde{A}_1\setminus \tilde{P}_1=\emptyset$, it is natural to set $\mathcal A_{\sigma_1}=\emptyset$ and $g_{\sigma_1}=e$, the identity element of $G$. For $j=2,$ we have $\tilde{A}_2\setminus \tilde{P}_2=\tilde{A}_1\cap \tilde{A}_2\subseteq \tilde{B}_1=\dot{g}_1\tilde{A}_1.$ Hence $\tilde{A}_2\setminus \tilde{P}_2=\dot{g}_1(\dot{g}_1^{-1}\tilde D_1^2)$. Set $\sigma_2=1,$ $\Lambda_2=\{\sigma_2\} $ and $g_{\sigma_2}=\dot{g}_1.$ Note that $|\Lambda_2|=1=2^{2-1}-1.$
Now let for $1\leq t\leq j$
$$ \tilde{A}_t\setminus \tilde{P}_t=\bigsqcup_{\sigma_t\in \Lambda_t}g_{\sigma_t}{\mathcal A}_{\sigma_t}\ \ \ \ \ \text{and}\ \ \ \ |\Lambda_t|=2^{t-1}-1.$$
By (\ref{b2})
\begin{align*}
  \tilde{A}_{j+1}\setminus \tilde{P}_{j+1} =& \bigsqcup_{t=1}^j\tilde D_t^{j+1} \\ =& \bigsqcup_{t=1}^j\dot{g}_t(\dot{g}_t^{-1}\tilde D_t^{j+1}) \\
   =& \bigsqcup_{t=1}^j\dot{g}_t\left([\dot{g}_t^{-1}\tilde D_t^{j+1}\cap \tilde{P}_{t}]\bigsqcup [\dot{g}_t^{-1}\tilde D_t^{j+1}\cap (\tilde{A}_{t}\setminus \tilde{P}_{t})] \right)\\
=& \bigsqcup_{t=1}^j\dot{g}_t\left([\dot{g}_t^{-1}\tilde D_t^{j+1}\cap \tilde{P}_{t}]\bigsqcup [\dot{g}_t^{-1}\tilde D_t^{j+1}\cap (\bigsqcup_{\sigma_t\in \Lambda_t}g_{\sigma_t}{\mathcal A}_{\sigma_t})] \right)\\
    = &\left( \bigsqcup_{t=1}^j\dot{g}_t(\dot{g}_t^{-1}\tilde D_t^{j+1}\cap \tilde{P}_{t}) \right) \bigsqcup \left(  \bigsqcup_{t=2}^j\bigsqcup_{\sigma_t\in \Lambda_t}\dot{g}_tg_{\sigma_t}[ g_{\sigma_t}^{-1}\dot{g}_t^{-1}\tilde D_t^{j+1}\cap {\mathcal A}_{\sigma_t}]\right).
\end{align*}
We have used the fact that $D_t^{j+1}\subseteq P_t$. Recall for the last equation that ${\mathcal A}_{\sigma_1}=\emptyset$. Note that the number of pieces in the decomposition above is $j+\sum_{t=2}^{j}|\Lambda_t|,$
so by the induction hypothesis $$|\Lambda_{j+1}|=j+\sum_{t=2}^{j}(2^{t-1}-1)=2^j-1.$$
Hence it is done for $j+1$.
 This completes the Step 2.
Before proceeding to the rest of the proof, we draw the reader's attention to the next two remarks.
\begin{remark}\label{remark}
In step 2 we explained how the sets ${\mathcal A}_{\sigma_{j+1}}$ are obtained from the sets $P_t$ and ${\mathcal A}_{\sigma_t},$ where $1\leq t<j+1$.
In fact each ${\mathcal A}_{\sigma_{j+1}}$ is in one of the following two types:
$$ \dot{g}_t^{-1}\tilde D_t^{j+1}\cap \tilde{P}_{t}\ \ \ \ \ \ \ \ \ \ \ \ \text{(I)}$$
$$  g_{\sigma_t}^{-1}\dot{g}_t^{-1}\tilde D_t^{j+1}\cap {\mathcal A}_{\sigma_t},\ \ \ \ \ \text{(II)}$$
 for suitable $t<j+1$ and $\sigma_t\in \Lambda_t$. So, one can easily see that ${\mathcal A}_{\sigma_{j+1}}\cap {\mathcal A}_{\sigma_t}\neq \emptyset$  if and only if  ${\mathcal A}_{\sigma_{j+1}}\subseteq {\mathcal A}_{\sigma_t}$. In this case for every $\sigma'_t\neq \sigma_t,$  ${\mathcal A}_{\sigma_{j+1}}\cap {\mathcal A}_{\sigma'_t}= \emptyset$.
\end{remark}
\begin{remark}\label{remark2}
For  $1\leq r<t\leq p$ and $\sigma_t\in \Lambda_t,$ either ${\mathcal A}_{\sigma_t}$ and $\tilde P_r$ are disjoint or ${\mathcal A}_{\sigma_t}\subseteq \tilde P_r$. To see this we apply Remark \ref{remark}. Let ${\mathcal A}_{\sigma_t}\cap \tilde P_r\neq \emptyset$. If ${\mathcal A}_{\sigma_t}$ is of type (I), then ${\mathcal A}_{\sigma_t}\subseteq \tilde P_k,$ for some $k<t.$  In this case $\tilde P_r\cap \tilde P_k\neq \emptyset,$ which is impossible unless $k=r$. Therefore ${\mathcal A}_{\sigma_t}\subseteq \tilde P_r$.
Now suppose that ${\mathcal A}_{\sigma_t}$ is of type (II) and $l$ be the smallest natural number such that for some $\sigma_l\in\Lambda_l$, ${\mathcal A}_{\sigma_t}\subseteq {\mathcal A}_{\sigma_l}$ and ${\mathcal A}_{\sigma_l}$ is of type (II). Then ${\mathcal A}_{\sigma_l}=\tilde P_m\cap \dot{g}^{-1}\tilde D_m^l\subseteq \tilde P_m$, for some $m<l$, which is not possible unless $m=r$. This implies that ${\mathcal A}_{\sigma_t}\subseteq \tilde P_r$.

\end{remark}

\noindent{\bf{Step 3.}} We show that for $1\leq t\leq p$  and $C\in Q_t$ there exists a family   $\{\mathcal M^C_{\delta_t},\ \ \ \delta_t\in \Gamma_t\}$ with cardinality $2^{t-1}$ of disjoint subsets of $\bigsqcup_{j=1}^t\tilde{P}_j$ and the members  $\{h^C_{\delta_t},\ \ \ \delta_t\in \Gamma_t\}$ of $G$ such that
\begin{equation}\label{bs}
x_0(C)=\bigsqcup_{\delta_t\in \Gamma_t}h^C_{\delta_t}\mathcal \mathcal M^C_{\delta_t}.\end{equation}
Since $Q_t\subseteq B_t,$ by Step 1 we have
\begin{align*}
  \tilde{Q}_t &\subseteq \tilde{B}_t=\dot{g}_t\tilde{A}_t \\
   & = \dot{g}_t\tilde{P}_t\bigsqcup \dot{g}_t(\tilde{A}_t\setminus \tilde P_t)\\
      & = \dot{g}_t\tilde{P}_t\bigsqcup \dot{g}_t(\bigsqcup_{\sigma_t\in \Lambda_t}g_{\sigma_t}{\mathcal A}_{\sigma_t})\\
        & = \dot{g}_t\tilde{P}_t\bigsqcup (\bigsqcup_{\sigma_t\in \Lambda_t}\dot{g}_tg_{\sigma_t}{\mathcal A}_{\sigma_t}).
\end{align*}
Clearly by Step 1 the number of pieces in this decomposition is $1+|\Lambda_t|=2^{t-1}$. Therefore for $C\in Q_t,$ one has
\begin{equation*}
  x_0(C)= \dot{g}_t(\tilde{P}_t\cap \dot{g}_t^{-1}x_0(C))\bigsqcup (\bigsqcup_{\sigma_t\in \Lambda_t}\dot{g}_tg_{\sigma_t}({\mathcal A}_{\sigma_t}\cap g_{\sigma_t}^{-1}\dot{g}_t^{-1}x_0(C))).
\end{equation*}
Now set $M^C_{\delta_0}=\tilde{P}_t\cap \dot{g}_t^{-1}x_0(C)$  and $M^C_{\delta_t}={\mathcal A}_{\sigma_t}\cap g_{\sigma_t}^{-1}\dot{g}_t^{-1}x_0(C))$, where $\sigma_t\in \Lambda_t$. Note that these sets are pairwise disjoint; $M^C_{\delta_0}$ is a subset of  $\tilde{P}_t$  whereas  the sets  $M^C_{\delta_t},\ \delta_t\in \Lambda_t$   are disjoint subsets of  $\bigsqcup_{j=1}^{t-1}\tilde{P}_j$  and besides   $(\bigsqcup_{j=1}^{t-1}\tilde{P}_j)\cap \tilde{P}_t=\emptyset$. Define $\Gamma_t=\{\delta_0\}\cup \Lambda_t$, $h_{\delta_t}^C=\dot{g}_tg_{\sigma_t},$ where $\delta_t=\sigma_t\in\Lambda_t$ and $h_{\delta_{t_0}}^C=\dot g_t$.  Thus $\{M^C_{\delta_t}, \delta_t\in\Gamma_t \}$   is a family of pairwise disjoint subsets of $G$ and
$$x_0(C)=\bigsqcup_{\delta_t\in \Gamma_t}h^C_{\delta_t}\mathcal \mathcal M^C_{\delta_t}.$$\\
It is noticeable that for $1\leq t\neq s\leq p,$ $M^C_{\delta_t}\cap M^{C'}_{\delta_s}=\emptyset$, where  $C\in Q_t$ and $C'\in Q_s$. We show this below. Without loss of generality assume that $t<s.$ Letting   $C\in Q_t$ and $C'\in Q_s$, we need to prove that the following four equalities are satisfied
\begin{align*}
  & M^C_{\delta_0}\cap M^{C'}_{\delta_0}=\emptyset, \\
   & M^C_{\delta_t}\cap M^{C'}_{\delta_0}=\emptyset, \\
   & M^C_{\delta_0}\cap M^{C'}_{\delta_s}=\emptyset, \\
   & M^C_{\delta_t}\cap M^{C'}_{\delta_s}=\emptyset.
\end{align*}
But it is equivalent to show that
\begin{align}
  & [\tilde{P}_t\cap \dot{g}_t^{-1}\tilde Q_t]\cap [\tilde{P}_s\cap \dot{g}_s^{-1}\tilde Q_s]=\emptyset, \\
   &[\tilde{P}_s\cap \dot{g}_s^{-1}\tilde Q_s]\cap [{\mathcal A}_{\sigma_t}\cap g_{\sigma_t}^{-1}\dot{g}_t^{-1}\tilde Q_t)]=\emptyset, \\
   & [\tilde{P}_t\cap \dot{g}_t^{-1}\tilde Q_t]\cap [{\mathcal A}_{\sigma_s}\cap g_{\sigma_s}^{-1}\dot{g}_s^{-1}\tilde Q_s)]=\emptyset, \label{3} \\
   & [{\mathcal A}_{\sigma_t}\cap g_{\sigma_t}^{-1}\dot{g}_t^{-1}\tilde Q_t)]\cap [{\mathcal A}_{\sigma_s}\cap g_{\sigma_s}^{-1}\dot{g}_s^{-1}\tilde Q_s)]=\emptyset.\label{4}
\end{align}
The first two equalities are trivial, since for each $r\neq s$, $P_r\cap P_s=\emptyset$ and besides, ${\mathcal A}_{\sigma_t}\subseteq \bigsqcup_{r=1}^{t-1} \tilde P_r$ and $t<s$.
We now prove (\ref{3}).
If ${\mathcal A}_{\sigma_s}$ is of type (I), then there exists $ l<s$ such that ${\mathcal A}_s= \dot{g}_l^{-1}\tilde D_l^s\cap \tilde P_l$. So
$$[\tilde{P}_t\cap \dot{g}_t^{-1}\tilde Q_t]\cap [{\mathcal A}_{\sigma_s}\cap g_{\sigma_s}^{-1}\dot{g}_s^{-1}\tilde Q_s)]\subseteq \tilde P_l\cap \tilde P_t.$$
The set $\tilde P_l\cap \tilde P_t$ is clearly empty if $l\neq t$. If $l= t,$
$$[\tilde{P}_t\cap \dot{g}_t^{-1}\tilde Q_t]\cap [{\mathcal A}_{\sigma_s}\cap g_{\sigma_s}^{-1}\dot{g}_s^{-1}\tilde Q_s]\subseteq \tilde P_l\cap \dot{g}_l^{-1}\tilde Q_l\cap \dot{g}_l^{-1}\tilde D_l^s$$ and $\tilde P_l\cap \dot{g}_l^{-1}\tilde Q_l\cap \dot{g}_l^{-1}\tilde D_l^s$ is also empty since $\tilde Q_l\cap \tilde D_l^s$, by (\ref{q_t}).
Now let ${\mathcal A}_{\sigma_s}$ be of type (II) and $k$ be the smallest natural number such that ${\mathcal A}_{\sigma_s}\subseteq {\mathcal A}_{\sigma_k}$ and ${\mathcal A}_{\sigma_k}$ is of type (II). This implies that ${\mathcal A}_{\sigma_k}\subseteq \tilde P_l\cap \dot{g}_l^{-1}\tilde D_l^k$, for some $l<k$. If (\ref{3}) does not satisfy, then by Remark \ref{remark2}, ${\mathcal A}_{\sigma_s}\subseteq \tilde P_t$
 and this is impossible unless $t=l$.  So ${\mathcal A}_{\sigma_s}\subseteq {\mathcal A}_{\sigma_k}\subseteq \dot{g}_t^{-1}\tilde D_t^k\cap \tilde P_t$. Taking into account that $D_t^k\cap \tilde Q_t=\emptyset$, we have
$$[\tilde{P}_t\cap \dot{g}_t^{-1}\tilde Q_t]\cap {\mathcal A}_{\sigma_s}\subseteq  \dot{g}_t^{-1}\tilde Q_t\cap \dot{g}_t^{-1}\tilde D_t^k=\emptyset.$$ This contradicts our assumption.\\
To prove (\ref{4}) let ${\mathcal A}_{\sigma_s}\cap {\mathcal A}_{\sigma_t}\neq \emptyset.$ By Remark \ref{remark}  ${\mathcal A}_{\sigma_s}\subseteq {\mathcal A}_{\sigma_t}$. In fact ${\mathcal A}_{\sigma_s}={\mathcal A}_{\sigma_t}\cap g_{\sigma_t}^{-1}\dot{g}_t^{-1}\tilde D_t^s.$
Thus we have
$$[{\mathcal A}_{\sigma_t}\cap g_{\sigma_t}^{-1}\dot{g}_t^{-1}\tilde Q_t]\cap [{\mathcal A}_{\sigma_s}\cap g_{\sigma_s}^{-1}\dot{g}_s^{-1}\tilde Q_s]\subseteq {\mathcal A}_{\sigma_t}\cap g_{\sigma_t}^{-1}\dot{g}_t^{-1}\tilde Q_t\cap g_{\sigma_t}^{-1}\dot{g}_t^{-1}\tilde D_t^s=\emptyset,$$ because $\tilde Q_t\cap \tilde D_t^s=\emptyset$.\\

\noindent{\bf{Step 4.}} $P_1,P_2,\dots,P_p$ are pairwise disjoint and for each
$C\in Con(\frak g,\mathcal E)$
$$\sum_{t=1}^p V'_t(C)\in \{0,1\}.$$
For $C\in Con(\frak g,\mathcal E)$ we put
$z_C=\alpha_C+\sum_{t=1}^p V'_t(C).$ This way, if $\sum_{t=1}^p V'_t(C)=1$, then $z_c\geq 2$. So, by (\ref{subtract}) there are at least two numbers $1\leq t^C_1,t^C_2\leq p$ such that $C\in Q_{t_1^C}$ and $C\in Q_{t_2^C}$. Therefore by Step 2,
$$x_0(C)=\bigsqcup_{\delta_{t_1^C}\in \Gamma_{t_1^C}}h^C_{\delta_{t_1^C}}\mathcal \mathcal M^C_{\delta_{t_1^C}}=\bigsqcup_{\delta_{t_2^C}\in \Gamma_{t_2^C}}h^C_{\delta_{t_2^C}}\mathcal \mathcal M^C_{\delta_{t_2^C}}.$$
On the other hand, if $\sum_{t=1}^p V'_t(C)=0,$ then $z_c\geq 1,$
so by (\ref{subtract}) there is at least one number $1\leq t_C\leq p$ such that $C\in Q_{t^C}$. Therefore by Step 2,
$$x_0(C)=\bigsqcup_{\delta_{t^C}\in \Gamma_{t^C}}h^C_{\delta_{t^C}}\mathcal \mathcal M^C_{\delta_{t^C}}.$$

Thus
\begin{multline*}
G\supseteq \left[ \bigsqcup_{\sum V'_i(C)=1}(\bigsqcup_{\delta_{t_1^C}\in \Gamma_{t_1^C}} \mathcal M^C_{\delta_{t_1^C}})\bigsqcup (\bigsqcup_{\delta_{t_2^C}\in \Gamma_{t_2^C}}\mathcal M^C_{\delta_{t_2^C}}) \right]\\
\bigsqcup \left[\bigsqcup_{\sum V'_i(C)=0}\bigsqcup_{\delta_{t^C}\in \Gamma_{t^C}}\mathcal M^C_{\delta_{t^C}}\right] \bigsqcup \left[\bigsqcup_{\sum V'_i(C)=0}x_0(C)\right].\end{multline*}
Observe that the sets in the right hand side of above inclusion are pairwise disjoint because
\begin{equation*}
 \{C\in \c,\ \sum V'_i(C)=0\}\subseteq\c\setminus\bigsqcup_{i=1}^p P_i,
\end{equation*}
i.e.,
\begin{equation*}
  \left(\bigsqcup_{\sum V'_i(C)=0}x_0(C) \right)\bigcap \left(\bigsqcup_{i=1}^p \tilde P_i \right)=\emptyset
\end{equation*}
and the sets
 $M^C_{\delta_t}, \delta_t\in\Gamma_t$ are pairwise disjoint subsets of $\bigsqcup_{i=1}^p \tilde P_i$ (take into account that $t_1^C,$ and $t_2^C$   are distinct numbers).
 Based on our choice of sets
\begin{align*}
G&= \left[ \bigsqcup_{\sum V'_i(C)=1}(\bigsqcup_{\delta_{t_1^C}\in \Gamma_{t^C}}h^C_{\delta_{t_1^C}}\mathcal \mathcal M^C_{\delta_{t_1^C}}) \right]\bigsqcup \left[\bigsqcup_{\sum V'_i(C)=0}(\bigsqcup_{\delta_{t^C}\in \Gamma_{t^C}}h^C_{\delta_{t^C}}\mathcal \mathcal M^C_{\delta_{t^C}}) \right] \\
G&=\left[ \bigsqcup_{\sum V'_i(C)=1}(\bigsqcup_{\delta_{t_2^C}\in \Gamma_{t_2^C}}h^C_{\delta_{t_2^C}}\mathcal \mathcal M^C_{\delta_{t_2^C}}) \right] \bigsqcup  \left[ \bigsqcup_{\sum V'_i(C)=0}x_0(C)\right],\end{align*}
which is a paradoxical decomposition of $G$ after omitting the empty sets from this decomposition.
This decomposition is complete if for each $C\in Con(\frak
g,\mathcal E)$, $$\sum_{i=1}^n V'_i(C)=a_C=1.$$ Otherwise
the decomposition is not complete but by a process described in
the proof of \cite[Proposition 1.2]{y-r} it can be changed into a
complete one.
\end{proof}

\begin{corollary}\label{tnumber}
 Using the notations of the proof of Theorem \ref{main},
if for every $C\in Con(\frak g,\mathcal E)$, $a_C=1,$ then
$\tau(G)\leq (\ell-1)(2^p-1).$
\end{corollary}
\begin{proof}
By the definition of configuration the cardinality of each $B_t$ cannot be more than $\ell-1$, where $|\c|=\ell$. The Step 3 in the proof of Theorem \ref{main} and the explanation following Proposition \ref{belg} yields
\begin{align*}
\tau(G)&\leq
\sum_{t=1}^{p}|\Lambda_t||Q_t|\\
&\leq   \sum_{t=1}^{p}2^{t-1}|B_t|\\
&\leq \sum_{t=1}^p2^{t-1}(\ell-1)\\
&\leq (\ell-1)(2^p-1).
\end{align*}
\end{proof}

\begin{remark}
In \cite{y-r} the authors constructed the paradoxical decomposition under the paradoxical condition. That is, a subsystem of $Eq(\frak g,\mathcal E)$ is equivalent to $BX=0,$ where each row of $B$ has nonnegative entries and  is of the form $\sum_{i=1}^m(L_i^{j_i}-L_i^{k_i})$, for some $j_i,k_i\in\{0,\dots,n\},$ where $L_i^j$ is the coefficient vector of the equation $$\sum_{x_j(C)\subseteq E_i}f_C-\sum_{x_0(C)\subseteq E_i}f_C=0$$
and in addition $B$ has no zero column.
Let $R^s=\sum_{i=1}^m(L_i^{j_i^s}-L_i^{k_i^s})$  be s-th  row of $B$. Then
$$\sum_{i=1}^{m}(\sum_{x_{j_i^s}(C)\subseteq E_i}f_C-\sum_{x_{k_i^s}(C)\subseteq E_i}f_C)=\sum_{C}R^s(C)f_C,$$ so
$$\sum_{i=1}^{m}(\sum_{x_{j_i^s}(C)\subseteq E_i}f_C)=\sum_{i=1}^{m}(\sum_{x_{k_i^s}(C)\subseteq E_i}f_C)+\sum_{C}R^s(C)f_C$$ therefore we have
$$\sum_{i=1}^{m}(\sum_{x_{0}(C)\subseteq E_i}f_C)=\sum_{i=1}^{m}(\sum_{x_{0}(C)\subseteq E_i}f_C)+\sum_{C}R^s(C)f_C.$$
But $\c=\bigsqcup_{i=1}^m \bigsqcup_{x_0(C)\subseteq E_i} x_0(C)$. Consequently we can replace $B$ by $W-V,$ where
$$ W=(\left( \begin{array}{cccc}
1,1,\dots,1 \\
1,1,\dots,1 \\
\vdots \\
1,1,\dots,1
\end{array} \right)+\left( \begin{array}{cccc}
R_1 \\
R_2 \\
\vdots \\
R_p
\end{array} \right))
\ \text{and}\ V=\left( \begin{array}{cccc}
1,1,\dots,1 \\
1,1,\dots,1 \\
\vdots \\
1,1,\dots,1
\end{array} \right).$$
Since by assumption each $R_s$ has nonnegative entries, $W_s-V_s-V_{s+1}$ has entries $\geq -1$. Finally $B$ has no zero column, so if $\sum_s({W_s-V_s})=(\alpha_C)_C,$ then $\alpha_C>0$.  This shows that the new system satisfies the normality condition in the statement of Theorem \ref{main}.
\end{remark}

\section{Graph interpretation}
In the current section we assign a graph to a group which helps us
construct the paradoxical decompositions and compute the Tarski
numbers.

\begin{definition}\label{graphic} \rm Let $G$ be a group and \cp $\Gamma=\Gamma(G,\frak g,\mathcal
E)$ is a graph constructed as follows\\ \\
$\centerdot$ The vertex
set of $\Gamma$ is identified with $\ell$-tuples of nonnegative integers, where $|\c|=\ell$:
 $$V(\Gamma):=\left\{(a_C)_{C\in Con(\frak g,\mathcal E)};\ \  a_C \in \mathbb{N}\cup\{0\}\right\}.$$\\
$\centerdot$
  There exists a directed edge from the vertex $A=(a_C)$ to
 $B=(b_C)$ only if $a_C\in\{0,1\}$ and there are disjoint subfamilies  $$\mathcal A^C_{1},\dots,\mathcal A^C_{{{b_C}}}\ \ \ \ (C\in \c)$$ of a partition of  $\bigsqcup_{a_C\neq 0}x_0(C)$ such that
 for all $C\in \c$
  $$x_0(C)=\bigsqcup_{A\in \mathcal A^C_{j}} g_A A,\ \ \ 1\leq j\leq b_C,$$
for suitable subsets $\{g_A,\ A\in \mathcal A^C_{j}\},$  $1\leq j\leq b_C,$ of $G$.
\end{definition}
\begin{proposition}\label{graph}
If $\Gamma$ contains adjacent vertices $A=(a_C)_{C\in \c}$ and  $B=(b_C)_{C\in \c}$  with
$a_C\in \{0,1\}$ and $\alpha_C:=b_C-a_C>0$, then $G$ admits a
paradoxical decomposition in terms of $\frak g$ and $\mathcal E$.
\end{proposition}
\begin{proof}
By the assumption for each $C$, $b_C\geq a_C+1$. In other words
$$b_C\geq\begin{cases}
   1, & \mbox{if } a_C=0 \\
    2, & \mbox{if } a_C\neq0.
  \end{cases}$$
Hence by Definition \ref{graphic} there exists a partition $\mathcal P$ of $\bigsqcup_{a_C\neq 0}x_0(C)$
with the following properties:
\begin{itemize}
  \item
If $a_C\neq 0$, there are $\mathcal A_1^C, \mathcal A_2^C\subseteq \mathcal P$  and $\{g_A,\ A\in \mathcal A_1^C\}, \{g_B,\ B\in \mathcal A_2^C\}\subseteq G$ such that
$$x_0(C)=\bigsqcup_{A\in \mathcal A_1^C}g_A.A=\bigsqcup_{B\in \mathcal A_2^C}g_B.B.$$

\item If $a_C= 0$, there are $\mathcal A^C\subseteq \mathcal P$ and $\{g_D,\ D\in \mathcal A^C\}\subseteq G$ such that
 $$x_0(C)=\bigsqcup_{D\in \mathcal A^C}g_D.D.$$
 \item All the above subfamilies of $\mathcal P$ are pairwise disjoint.
\end{itemize}
 For convenience denote the set $\bigsqcup_{a_C= 0}x_0(C)$ by $E$. We know that $G=\bigsqcup_{C\in \c}x_0(C)$. Then we have
$$G=E\bigsqcup \left( \bigsqcup_{a_C\neq 0}\bigsqcup_{A\in \mathcal A_1^C}g_A.A\right)$$ and
$$G= \left(\bigsqcup_{a_C= 0}\bigsqcup_{D\in \mathcal A^C}g_D.D\right)\bigsqcup \left( \bigsqcup_{a_C\neq 0}\bigsqcup_{B\in \mathcal A_2^C}g_B.B\right)$$
which give a (not necessarily complete) paradoxical decomposition for $G$. To make it complete, one can apply the proof of \cite[Proposition 1.2]{y-r} as usual.
 \end{proof}
\begin{theorem}
Let $G$ be a group and \cp If $Eq(\frak g,\mathcal E)$ has no
nonnegative nonzero solution with a normal subsystem, then $\Gamma=\Gamma(G,\frak g,\mathcal
E)$ includes vertices $A=(a_C)_{C\in \c}$ and  $B=(b_C)_{C\in \c}$  with
$a_C\in \{0,1\}$ and $\alpha_C:=b_C-a_C>0$.
\end{theorem}
\begin{proof}
  Using the notations
of the proof of Theorem \ref{main},  $A=\sum V'_i$ and $B=\sum
W'_i,$ are desired vertices because first, $B-A>0$ and $\sum V'(C)\in \{0,1\}$, and second, by (\ref{bs}) they satisfy the second condition of the definition of $\Gamma(G, \frak g,\mathcal E)$.
\end{proof}

\subsection{Diagrams associated with configuration equations}
Here we present a kind of diagram associated with configuration equations. It could be of great importance in order to reduce the complexity of the proofs in the previous section.
\begin{itemize}
  \item Let  $A, B_1,\dots,B_t$   be subsets of $\c$.
 If there is  $g\in G$ and there is a partition $\{A_1,\dots,A_t\}$ for $\tilde{A}=\bigsqcup_{C\in A}x_0(C)$  such that
  $$\tilde{B}_j=g{A_j},\ \ \ 1\leq j\leq t,$$
we use the following tree diagram
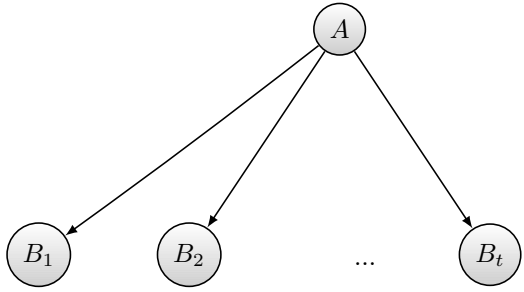
{{
\begin{center}
\begin{tikzpicture}[-latex ,auto ,node distance =3 cm and 2cm ,on grid ,
semithick , state/.style ={ circle ,top color =white , bottom
color = black!20 , draw,black , text=black, minimum width =.3cm}]

\node[state] (A) {$A$};
\node[state] (b_2) [below left=of A] {$B_2$};
\node[state] (b_1) [left  =of b_2] {$B_1$};
\node[state] (b_t) [below right =of A] {$B_t$};

\path (A) edge [bend left =1] node[below =0.15 cm] {$$} (b_1);

\path (A) edge[bend right =0] node[below =0.15 cm] {$$} (b_2);

\path (A) edge[bend right =0] node[below =1.5 cm] {$\ \ \ \ ...\ \ \ \ \ \ \ \ \ \ \ \ \ \ \  $} (b_t);.
\end{tikzpicture}
\end{center}
}}
\noindent{Clearly} if $D\subseteq A,$ and $E_j=\tilde B_j\cap g\tilde D$, we also have a similar diagram with directed paths from $D$ to $E_1,\dots ,E_{t-1},E_t$.
  \item By the diagram \\
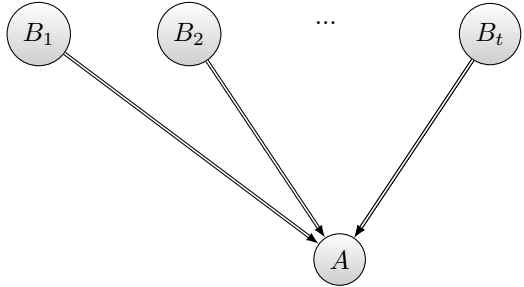
\begin{center}
\begin{tikzpicture}[-latex ,node distance =3 cm and 2cm ,on grid ,
state/.style ={ circle ,top color =white , bottom
color = black!20 , draw,black , text=black, minimum width =.3cm}]

\node[state] (A) {$A$};
\node[state] (b_2) [above left=of A] {$B_2$};
\node[state] (b_1) [left  =of b_2] {$B_1$};
\node[state] (b_t) [above right =of A] {$B_t$};
\path (b_2)  node[below =-.3 cm] {$  \ \ \ \ \ \  \ \ \ \ \ \ \ \ \ \ \ \ \ \ \ \  \ \ \ \ \ \ \ \ \ \ ...\  $} (b_t);.
\begin{scope}[edge style/.style={draw=black,double=white}]
 \draw[edge style] (b_t)--(A);
       \draw[edge style] (b_1)--(A);
       \draw[edge style] (b_2)--(A);
               \draw[edge style] (b_t)--(A);
               \end{scope}
\end{tikzpicture}
\end{center}
we mean $A=\bigsqcup_{i=1}^tB_i$.
  \item
Let $A,B\subseteq \c$. In the sequel we denote a directed path of any kind of above two types from $A$ to $B$ by $[A,B]$. Inductively $[A,B_1,B_2,\dots,B_n]$ is a path which connects the sets $A, B_1,B_2,\dots,B_n$ through the sequence of $n$ single paths  $[A,B_1],[B_1,B_2],\dots,[B_{n-1},B_n]$.
 \end{itemize}
In the proof of Theorem {\ref{main}} according to the normality condition we have the possibility to borrow the configurations in $ A_t\setminus P_t$ from $B_1,\dots,B_{t-1}$  so that the (not necessarily disjoint) subsets $  A_1,\dots,  A_p$ turn into disjoint subsets $  P_1,\dots,  P_p$ and
it is enough for us to construct a paradoxical decomposition of $G$ and to get an upper bound for $\tau(G)$. We Apply the above diagrams to make the proof of that theorem more clear to the reader. Using the same notations, by (\ref{q_t}) we have
$\tilde P_1=\tilde A_1=\dot g_1^{-1} \tilde B_1=\dot g_1^{-1}\tilde Q_1\bigsqcup \dot g_1^{-1}\tilde D_1^2\bigsqcup\dots \bigsqcup\dot g_1^{-1}\tilde D_1^p.$ So, there are $p$ directed paths from $P_1=A_1$ to $D_1^2,\dots, D_1^p$ and $Q_1.$
Also
$\tilde  D_1^2\bigsqcup \tilde P_2=\tilde A_2=\dot g_2^{-1} \tilde B_2= \dot g_2^{-1}\tilde Q_2\bigsqcup \dot g_2^{-1}\tilde D_2^3\bigsqcup\dots \bigsqcup\dot g_2^{-1}\tilde D_2^p.$ So, there are $p-1$ directed paths from $D^2_1$ to $D_2^3,\dots, D_2^p$ and $Q_2$ and the same number of paths from $P_2$ to $D_2^3,\dots, D_2^p$ and $Q_2$.
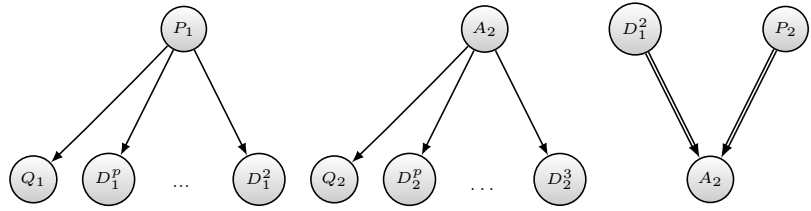
\begin{center}
\begin{tikzpicture}[-latex ,auto ,node distance =2 cm and 1cm ,on grid ,
semithick , state/.style ={ circle ,top color =white , bottom
color = black!20 , draw,black , text=black, minimum width =.6cm}]
\clip(-4,1) rectangle (9,-3);
{\tiny{
\node[state] (A) {$P_1$};
\node[state] (b_2) [below left=of A] {$D_1^p$};
\node[state] (b_1) [left  =of b_2] {$Q_1$};
\node[state] (b_t) [below right =of A] {$D_1^2$};
\node[inner sep=0,minimum size=0,right of=b_t]  (kk4) {}; 

\node[state] (b_0) [above right =of kk4] {$A_2$};

\path (A) edge [bend left =1] node[below =0.15 cm] {$$} (b_1);
\path (A) edge[bend right =0] node[below =0.15 cm] {$$} (b_2);
\path (A) edge[bend right =0] node[below =1 cm] {$\ ...\ \ \ \ \ \ \ \ \ \ \ \ \ $} (b_t);

\node[state] (b'_2) [below left=of b_0] {$D_2^p$};
\node[state] (b'_1) [left  =of b'_2] {$Q_2$};
\node[inner sep=0,minimum size=0,right of=b'_1]  (kkk4) {}; 
\node[state] (b'_t) [right =of kkk4] {$D_2^3$};

\path (b_0) edge [bend left =1] node[below =3.5 cm] {} (b'_1);
\path (b_0) edge[bend right =0] node[below =-3.1 cm] {} (b'_2);
 \path (b_0) edge[bend right =0] node[below =1 cm] {$\ \ \ \ \ \ \cdots\ \ \ \ \ \ \ \ \ \ \ \ \ \ \ \ \ \ $} (b'_t);.

\begin{scope} [edge style/.style={draw=black,double=white}]
\node[state] (bb_t) [above right=of b'_t] {$D_1^2$};
\node[state] (bbb_t) [below right=of bb_t] {$A_2$};
\node[state] (bbbb_t) [above right=of bbb_t] {$P_2$};

              \draw[edge style] (bbbb_t)--(bbb_t);
              \draw[edge style] (bb_t)--(bbb_t);
\end{scope}
}}
\end{tikzpicture}
\end{center}
Gluing these diagrams we obtain two paths from $D_1^2$ and $P_2$ to  $Q_2$ and so on.\vspace{.2cm}
{{
\def\radius{3.mm}
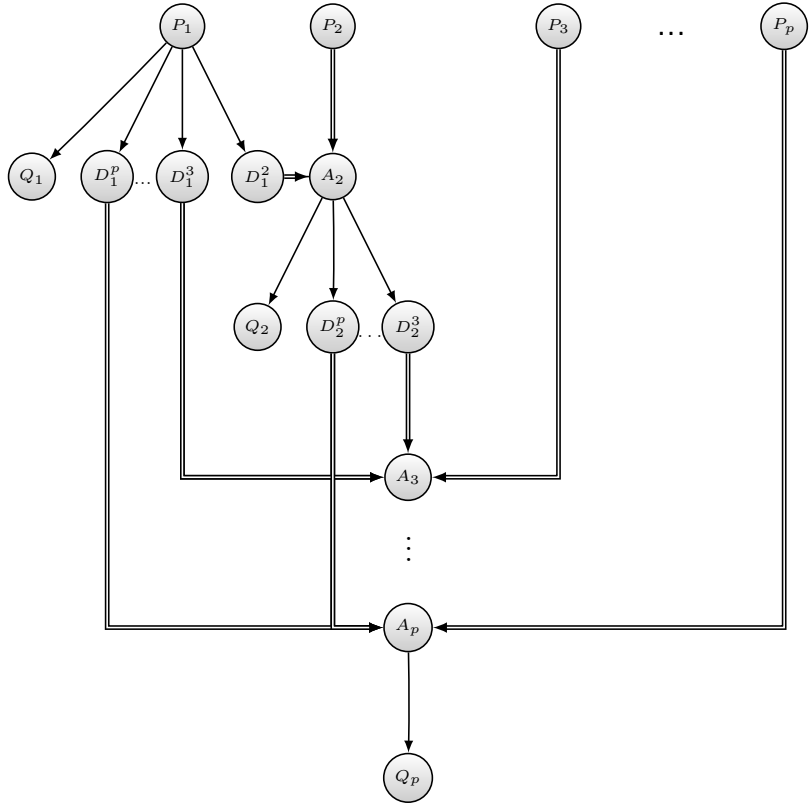
\begin{tikzpicture}[-latex ,auto ,node distance =2 cm and 1cm ,on grid ,
semithick , state/.style ={ circle ,top color =white , bottom
color = black!20 , draw,black , text=black, minimum width =.3cm}]
{\tiny{
\node[state] (A) {$P_1$};
\node[state] (b_2) [below left=of A] {$D_1^p$};
\node[state] (b_1) [left  =of b_2] {$Q_1$};
\node[state] (b_t) [below right =of A] {$D_1^2$};
\node[state] (b_0) [right=of b_t] {$A_2$};
\node[state] (p_0) [above=of b_0] {$P_2$};
\node[state] (pp_0) [left=of b_t] {$D_1^3$};

\path (A) edge [bend left =1] node[below =0.15 cm] {$$} (b_1);
{\draw[-] (b_t) edge node {$$} (b_0);}
\path (A) edge[bend right =0] node[below =0.15 cm] {$$} (b_2);
\path (A) edge[bend right =0] node[below =1 cm] {$\ ...\ \ \ \ \ \ \ \ \ \ \ \ \ $} (pp_0);
\path (A) edge[bend right =0] node[below =1.5 cm] {$$} (b_t);

\node[state] (b'_2) [below left=of b_0] {$Q_2$};
\node[state] (b'_1) [below  =of b_0] {$D_2^p$};
\node[state] (b'_t) [below right =of b_0] {$D_2^3$};

\node[inner sep=0,minimum size=0,right of=p_0] (l1) {}; 
\node[state] (pppp_0) [ right =of l1] {$P_3$};
\node[inner sep=0,minimum size=0,right of=pppp_0] (l2) {}; 

\node[state] (ppppp_0) [ right =of l2] {$P_p$};

\node[inner sep=0,minimum size=0,below of=pppp_0] (k1) {}; 
\node[inner sep=0,minimum size=0,below of=k1] (k2) {}; 
\node[inner sep=0,minimum size=0,below of=k2] (k3) {}; 

\node[inner sep=0,minimum size=0,below of=pp_0] (k5) {}; 
\node[inner sep=0,minimum size=0,below of=k5] (k6) {}; 
\node[inner sep=0,minimum size=0,below of=b'_1]  (kk4) {}; 

\node[state] (k4) [below =of b'_t] {$A_3$};
\node[state] (p4) [below =of k4] {$A_p$};
\node[state] (qp) [below =of p4] {$Q_p$};
\path (p4) edge [bend left =1] node[below =3.5 cm] {$$} (qp);

\path (b_0) edge [bend left =1] node[below =3.5 cm] {$\hspace{2cm}{\boldmath{\Large{\vdots}}} $} (b'_1);
\path (b_0) edge[bend right =0] node[below =-3.1 cm] {$\hspace{10cm}{\boldmath{\Huge{\centerdot\centerdot\centerdot}}} $} (b'_2);
\begin{scope} [edge style/.style={draw=black,double=white}]
              \draw[edge style] (p_0)--(b_0);
       \draw[edge style] (b_t)--(b_0);
       \draw[edge style] (pppp_0)|-(k4);
               \draw[edge style] (b'_t)--(k4);

        \draw[edge style, name path=line 1] (pp_0)|-(k4);

\draw[edge style, name path=line 2] (b'_1)|-(p4);

\path [name intersections={of = line 1 and line 2}];
  \coordinate (S)  at (intersection-1);

  \path[name path=circle] (S) circle(\radius);

  \path [name intersections={of = circle and line 2}];
  \coordinate (I1)  at (intersection-1);
  \coordinate (I2)  at (intersection-2);

 \draw[edge style, name path=line 1] (pp_0)|-(k4);

     \end{scope}

     \begin{scope} [edge style/.style={draw=black,double=white}]
             \draw[edge style] (ppppp_0)|-(p4);
       \draw[edge style] (b_2)|-(p4);
\draw[edge style, name path=line 2] (b'_1)|-(p4);

\path [name intersections={of = line 1 and line 2}];
  \coordinate (S)  at (intersection-1);

  \path[name path=circle] (S) circle(\radius);

  \path [name intersections={of = circle and line 2}];
  \coordinate (I1)  at (intersection-1);
  \coordinate (I2)  at (intersection-2);
%

     \end{scope}
   \path (b_0) edge[bend right =0] node[below =1 cm] {$\ \ \ \ \ \ \ \ \ \ \ \ \ \ \ \ \ \ \cdots\ \ \ \ \ \ \ \ \ \ \ \ \ \ \ \ \ \ $} (b'_t);.
 }}
\end{tikzpicture}

 Our diagram contains the following $2^p-1$ paths
\begin{eqnarray*}
&&  [P_1,Q_1],  \\
&&   [P_1,D_1^2,A_2,Q_2], [P_2,A_2,Q_2],  \\
    &&[P_1,D_1^3,A_3,Q_3],[P_1,D_1^2,A_2,D_2^3,A_3,Q_3],[P_2,A_2,D_2^3,A_3,Q_3], [P_3,A_3,Q_3],\\
    && \vdots  \\
   && [P_1,D_1^p,A_p,Q_p],\dots,[P_1,D_1^2,A_2,D_2^p,A_p,Q_p],\dots,[P_p,A_p,Q_p].
\end{eqnarray*}
Depending on the members of $Q_1,\dots,Q_p,$ these paths can be continued one more step.
In fact $\bigsqcup_{t=1}^p P_t$ admits a partition with at most $2^p-1$ pieces such that with suitable coefficients construct the sets $\bigsqcup_{C\in Q_t}x_0(C),\ 1\leq t\leq p$, separately (see proof of Theorem \ref{main}). $2^{t-1}$ pieces in the partition are used to make $\tilde Q_t.$ These pieces are divided into $2^{t-1}|Q_t|$ new pieces to construct all the sets $x_0(C),\ C\in Q_t$. To have a paradoxical decomposition each  $x_0(C),\ C\in \c$  must be obtained at least two times. If $\sum_{t=1}^p V'_t(C)=1,$ then $\sum_{t=1}^p W'_t(C)\geq 2.$ In other words, $C$ belongs to at least two $Q_ts$ and it is done. Otherwise, $C$ may appear in only one $Q_t$ and we have to add the single set ${\{C,\ \sum_{t=1}^p V'_t(C)=0\tilde\}}$ to the pieces of our paradoxical decomposition. Therefore $\tau(G)\leq \sum_{t=1}^p 2^{t-1}|Q_t|\bf{+1}$. It is to be noted that $\sum_{t=1}^p 2^{t-1}|Q_t|$ is the maximum possible number of mentioned paths.\\

 Naturally to make this observation more accurate, we have to choose a procedure in which the number of pieces in the paradoxical decomposition is as small as possible. This process depends directly on the set of configurations given in the statement of problem. We explain it in the next subsection.

\subsection{Examples}
We end the paper with some examples. Our examples will be based on the sets of configurations instead of the elements of the group.
As the paradoxical decomposition of
a group is not unique, different diagrams also exist.
One can find a more precise Tarski number's upper bound by counting the paths from the
top to the bottom in a minimal diagram (a diagram with minimum number of paths). It is noticeable that a minimal diagram is not necessarily a connected one.
\begin{example}\label{tt5}
Let $C_1=(1,2,3,2),\ C_2=(1,3,1,3),\ C_3=(2,1,2,2),\
C_4=(3,3,1,2)$ and $\ C_5=(3,3,2,1)$ be  configurations corresponding to the configuration pair $(\frak g,\mathcal E)$ of a
group $G$.
Consider the following equations
\begin{align*}
\sum_{x_1(C)\subseteq E_1}f_C =&\sum_{x_0(C)\subseteq E_1}f_C,\\   \sum_{x_1(C)\subseteq E_1}f_C=&\sum_{x_0(C)\subseteq E_1}f_C\\
\sum_{x_3(C)\subseteq E_1}f_C=&\sum_{x_0(C)\subseteq E_1}f_C,\\     \sum_{x_3(C)\subseteq E_1}f_C=&\sum_{x_0(C)\subseteq E_1}f_C\\
\sum_{x_0(C)\subseteq E_2}f_C=&\sum_{x_2(C)\subseteq E_2}f_C,\\     \sum_{x_1(C)\subseteq E_2}f_C=&\sum_{x_2(C)\subseteq E_2}f_C\\
\sum_{x_3(C)\subseteq E_3}f_C=&\sum_{x_0(C)\subseteq E_3}f_C
\end{align*}
which are the next equations, respectively
\begin{align*}
f_{C_3}=&f_{C_1}+f_{C_2},\\  f_{C_3}=&f_{C_1}+f_{C_2},\\ f_{C_5}=&f_{C_3},\\ f_{C_5}=&f_{C_3},\\ f_{C_3}=&f_{C_3}+f_{C_5},\\ f_{C_1}=&f_{C_3}+f_{C_5},\\ f_{C_2}=&f_{C_4}+f_{C_5}.
\end{align*}
Then  $A=\left( \begin{array}{ccccc}
0&0&1&0&0 \\
0&0&1&0&0 \\
0&0&0&0&1 \\
0&0&0&0&1 \\
0&0&1&0&0 \\
1&0&0&0&0 \\0&1&0&0&0
\end{array} \right)\ \  \text{and}\ \ \  B=\left( \begin{array}{ccccc}
1&1&0&0&0 \\
1&1&0&0&0 \\
0&0&1&0&0 \\
0&0&1&0&0 \\
0&0&1&0&1 \\
0&0&1&0&1\\0&0&0&1&1
\end{array} \right)$ are the corresponding coefficient matrices.
Setting $\pi=(2\ 7)(4\ 6),$ one has
\begin{center}
$TP_{\pi}(B-A)-P_{\pi}^{+}A=\left(
\begin{array}{ccccc}
1&0&-1&0&0 \\
1&0&-1&1&0 \\
0&0&0&1&0 \\
0&0&0&1&1 \\
0&0&1&1&1 \\
0&0&1&1&1 \\1&1&0&1&1
\end{array} \right).$\\
\end{center}

Therefore this system is normal.
The equations corresponding to $TP_{\pi}(B-A)$ are
\begin{align*}
f_{C_3}=&f_{C_1}+f_{C_2},\\ f_{C_2}=&f_{C_4}+f_{C_5} ,\\ f_{C_5}=&f_{C_3},\\ f_{C_1}=&f_{C_3}+f_{C_5}, \\ f_{C_3}=&f_{C_3}+f_{C_5},\\ f_{C_5}=&f_{C_3},  \\ f_{C_3}=&f_{C_1}+f_{C_2}.
\end{align*}
We use the notations of the proof of Theorem \ref{main}. In this example $p=7$, $A_1=\{C_3\}, A_2=\{C_2\}, A_3=\{C_5\} ,A_4=\{C_1\}, A_5=\{C_3\}, A_6=\{C_5\}, A_7=\{C_3\}$, $B_1=\{C_1, C_2\}, B_2=\{C_4, C_5\}, B_3=\{C_3\}, B_4=\{C_3, C_5\}, B_5=\{C_3, C_5\}, B_6=\{C_3\}, B_7=\{C_1, C_2\}.$
The initial diagram is
 {\tiny
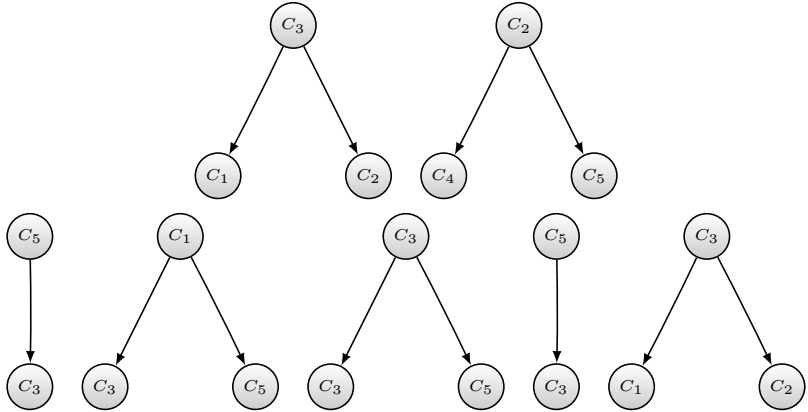
{
\begin{center}
\begin{tikzpicture}[-latex ,auto ,node distance =2 cm and 1cm ,on grid ,
semithick , state/.style ={ circle ,top color =white , bottom
color = black!20 , draw,black , text=black, minimum width =.3cm}]
\node[state] (Cc) {$C_5$};
\node[state] (Aa) [above left=of Cc] {$C_2$};
\node[state] (Bb) [below left =of Aa] {$C_4$};
\node[state] (C) [left=of Bb]{$C_2$};
\node[state] (A) [above left=of C] {$C_3$};
\node[state] (B) [below left=of A] {$C_1$};
\path (A) edge [bend left =1] node[below =0.15 cm] {$$} (B);

\path (A) edge[bend right =0] node[below =0.15 cm] {$$} (C);

\path (Aa) edge [bend left =1] node[below =0.15 cm] {$$} (Bb);
\path (Aa) edge[bend left =0] node[below =0.15 cm] {$$} (Cc);

\end{tikzpicture}
\end{center}
}}
{\tiny{
\begin{center}
\begin{tikzpicture}[-latex ,auto ,node distance =2 cm and 1cm ,on grid ,
semithick , state/.style ={ circle ,top color =white , bottom
color = black!20 , draw,black , text=black, minimum width =.3cm}]
\node[state] (Ccc) {$C_2$};
\node[state] (Aaa) [above left=of Ccc] {$C_3$};
\node[state] (Bbb) [below left =of Aaa] {$C_1$};
\node[state] (Cccc)  [left=of Bbb]{$C_3$};
\node[state] (Aaaa) [above =of Cccc] {$C_5$};
\node[state] (Ccccc)  [left=of Cccc]{$C_5$};
\node[state] (Aaaaa) [above left =of Ccccc] {$C_3$};
\node[state] (Bbbbb) [below left =of Aaaaa] {$C_3$};
\node[state] (Cccccc)  [left=of Bbbbb]{$C_5$};
\node[state] (Aaaaaa) [above left=of Cccccc] {$C_1$};
\node[state] (Bbbbbb) [below left =of Aaaaaa] {$C_3$};
\node[state] (Ccccccc)  [left=of Bbbbbb]{$C_3$};
\node[state] (Aaaaaaa) [above =of Ccccccc] {$C_5$};

\path (Aaa) edge [bend left =1] node[below =0.15 cm] {$$} (Bbb);
\path (Aaa) edge [bend left =1] node[below =0.15 cm] {$$} (Ccc);

\path (Aaaa) edge [bend left =1] node[below =0.15 cm] {$$} (Cccc);

\path (Aaaaa) edge[bend right =0] node[below =0.15 cm] {$$} (Bbbbb);

\path (Aaaaa) edge [bend left =1] node[below =0.15 cm] {$$} (Ccccc);

\path (Aaaaaa) edge[bend right =0] node[below =0.15 cm] {$$} (Bbbbbb);

\path (Aaaaaa) edge [bend left =1] node[below =0.15 cm] {$$} (Cccccc);

\path (Aaaaaaa) edge [bend left =1] node[below =0.15 cm] {$$} (Ccccccc);
\end{tikzpicture}
\end{center}
}}
Applying the process in  the proof of Theorem, we have  $P_1=\{C_3\}, P_2=\{C_2\}, P_3=\{C_5\} ,P_4=\{C_1\}, P_5=\emptyset, P_6=\emptyset, P_7=\emptyset$, $Q_1=\{C_1, C_2\}, Q_2=\{C_4\}, Q_3=\emptyset, Q_4=\{C_3, C_5\}, Q_5=\{ C_5\}, Q_6=\{C_3\}, Q_7=\{C_1, C_2\}.$ The associated diagram is

 {\tiny
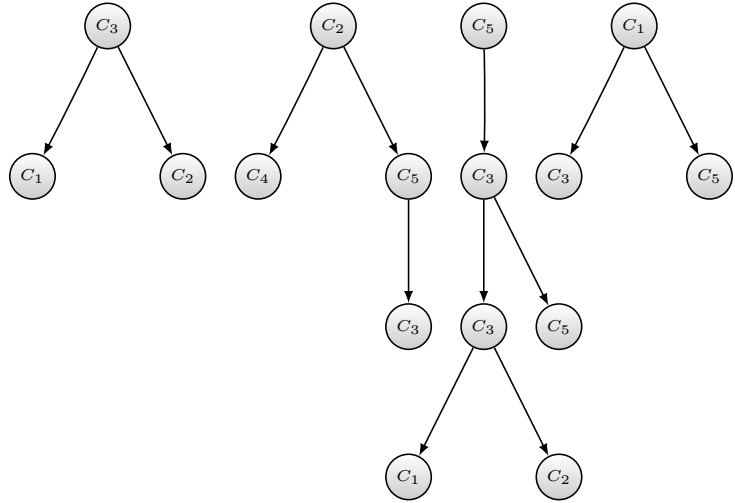
{
\begin{center}
\begin{tikzpicture}[-latex ,auto ,node distance =2 cm and 1cm ,on grid ,
semithick , state/.style ={ circle ,top color =white , bottom
color = black!20 , draw,black , text=black, minimum width =.3cm}]
\node[state] (Cc) {$C_5$};
\node[state] (Aa) [above left=of Cc] {$C_2$};
\node[state] (Bb) [below left =of Aa] {$C_4$};
\node[state] (C) [left=of Bb]{$C_2$};
\node[state] (A) [above left=of C] {$C_3$};
\node[state] (B) [below left=of A] {$C_1$};
\node[state] (Bbb) [right=of Cc] {$C_3$};
\node[state] (Aaa) [above =of Bbb] {$C_5$};
\node[state] (Bbbb) [right=of Bbb] {$C_3$};
\node[state] (Ccc) [above right=of Bbbb] {$C_1$};
\node[state] (Bbbbb) [below right =of Ccc] {$C_5$};
\node[state] (Bbbbbb) [below =of Cc] {$C_3$};
\node[state] (Bc) [below =of Bbb] {$C_3$};
\node[state] (Ac) [right =of Bc] {$C_5$};

\node[state] (Bcc) [below left=of Bc] {$C_1$};
\node[state] (Acc) [below right =of Bc] {$C_2$};

\path (A) edge [bend left =1] node[below =0.15 cm] {$$} (B);
\path (Aaa) edge [bend left =1] node[below =0.15 cm] {$$} (Bbb);

\path (A) edge[bend right =0] node[below =0.15 cm] {$$} (C);

\path (Aa) edge [bend left =1] node[below =0.15 cm] {$$} (Bb);
\path (Aa) edge[bend left =0] node[below =0.15 cm] {$$} (Cc);
\path (Ccc) edge [bend left =1] node[below =0.15 cm] {$$} (Bbbb);
\path (Ccc) edge[bend left =0] node[below =0.15 cm] {$$} (Bbbbb);
\path (Cc) edge[bend left =0] node[below =0.15 cm] {$$} (Bbbbbb);

\path (Bbb) edge[bend left =0] node[below =0.15 cm] {$$} (Bc);
\path (Bbb) edge[bend left =0] node[below =0.15 cm] {$$} (Ac);

\path (Bc) edge[bend left =0] node[below =0.15 cm] {$$} (Bcc);
\path (Bc) edge[bend left =0] node[below =0.15 cm] {$$} (Acc);
\end{tikzpicture}
\end{center}
}}
The condition of normality enables us to have disjoint subsets $P_1,\dots, P_p$ which produce $\c$ two times (see proof of Theorem \ref{main}). By Corollary
\ref{tnumber}, $\tau(G)\leq (\ell-1)(2^p-1)=508$ ($\ell=5$ and $p=7$). This is not a very accurate upper bound. Changing the diagram to a minimal one, helps us compute the most precise bound  $\tau(G)$. Look at the following diagram
{\tiny
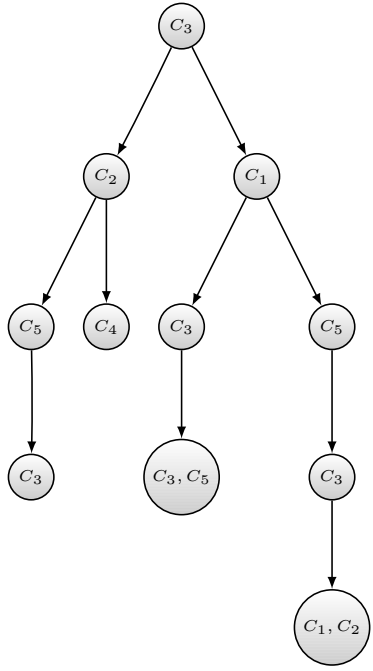
{
\begin{center}
\begin{tikzpicture}[-latex ,auto ,node distance =2 cm and 1cm ,on grid ,
semithick , state/.style ={ circle ,top color =white , bottom
color = black!20 , draw,black , text=black, minimum width =.3cm}]

\node[state] (C) {$C_1$};
\node[state] (A) [above left=of C] {$C_3$};
\node[state] (B) [below left =of A] {$C_2$};
\node[state] (D) [below right=of C] {$C_5$};
\node[state] (H) [below left=of C] {$C_3$};
\node[state] (M) [black,below =of B] {$C_4$};
\node[state] (N) [ left=of M] {$C_5$};
\node[state] (G) [black,below  =of H] {$\tiny{C_3,C_5}$};
\node[state] (O) [below  =of D] {$C_3$};
\node[state] (R) [black,below  =of O] {$C_1, C_2$};
\node[state] (Q) [below =of N] {$C_3$};

\path (A) edge [bend left =1] node[below =0.15 cm] {$$} (B);

\path (A) edge[bend right =0] node[below =0.15 cm] {$$} (C);

\path (C) edge[bend right =0] node[below =0.15 cm] {$$} (D);

\path (C) edge[bend right =0] node[below =0.15 cm] {$$} (H);

\path (H) edge[bend right =0] node[below =0.15 cm] {$$} (G);

\path (D) edge[bend right =0] node[below =0.15 cm] {$$} (O);

\path (B) edge [bend left =0] node[below =0.15 cm] {$$} (M);

\path (B) edge[bend right =0] node[below =0.15 cm] {$$} (N);

\path (N) edge [bend left =1] node[below =0.15 cm] {$$} (Q);

\path (O) edge[bend right =0] node[below =0.15 cm] {$$} (R);
\end{tikzpicture}
\end{center}
}}
\noindent According to (\ref{key}) we have
$g_1x_0(C_3)= x_0(C_1)\bigsqcup x_0(C_2), g_3x_0(C_2)= x_0(C_4)\bigsqcup x_0(C_5)$, $g_2^{-1} g_1x_0(C_1)= x_0(C_3)\bigsqcup x_0(C_5)$, $g_2^{-1}x_0(C_3)= x_0(C_3)\bigsqcup x_0(C_5)$ and  $g_3x_0(C_5)= x_0(C_3).$

Then $x_0(C_3)=\Delta_1\bigsqcup\Delta_2\bigsqcup\Delta_3\bigsqcup\Delta_4,$ where
\begin{align*}
 \Delta_1&=g_1^{-1}g_3^{-2}x_0(C_3)\\ \Delta_2&= g_1^{-1}g_3^{-1}x_0(C_4) \\ \Delta_3&=g_1^{-2}g_2^{2}\left(x_0(C_3)\bigcup x_0(C_5)\right)\\
  \Delta_4&= g_1^{-2}g_2g_3^{-1}g_1^{-1}\left(x_0(C_1)\bigcup x_0(C_2)\right).\end{align*}
Therefore it is obvious that
 \begin{align*}
   G &=x_0(C_3)\bigsqcup(G\setminus x_0(C_3)) \\
     &=g_3g_1\Delta_2\bigsqcup(G\setminus x_0(C_3))
   \end{align*}
   and besides
    $$G=g_3^2g_1\Delta_1\bigsqcup g_2^{-2}g_1^{2}\Delta_3\bigsqcup g_1g_3g_2^{-1}g_1^2\Delta_4.$$
In other words, we have a paradoxical decomposition consisting of 5 pieces. Thus $\tau(G)\leq 5.$  This number can be obtained from the above diagram by counting the paths from the beginning points to the end points and adding the result by 1.
The suitable part of $\Gamma(G,\frak g,\mathcal E)$ corresponding to this example is\\
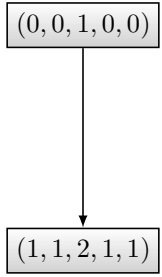
\begin{center}
\begin{tikzpicture}[-latex ,auto ,node distance =3 cm and 2cm ,on grid ,
semithick , state/.style ={ ,top color =white , bottom
color = black!20 , draw,black , text=black, minimum width =.3cm}]

\node[state] (C) {$(0,0,1,0,0)$};
\node[state] (A) [below =of C] {$(1,1,2,1,1)$};

\path (C) edge[bend right =0] node[below =0.15 cm] {$$} (A);
\end{tikzpicture}
\end{center}
\end{example}

We have a very special case of Theorem \ref{main} which
eventuates a result for computing the Tarski number. We use the notations of section
3.

\begin{theorem}
Suppose that $\sum_{i=1}^p(W_i-V_i)=(1,1,\dots,1)$ and there
exists a permutation matrix $P_{\pi}$ such that the first $p-1$ rows of
$P_{\pi}B-P_{\pi}^+A$ have nonnegative entries. Then $\tau(G)\leq
1+|A_{\pi(1)}|+\ell$.
\end{theorem}
\begin{proof}
The assumption says that for $1\leq i\leq n,$ $A_{\pi(i+1)}\subseteq
B_{\pi(i)}.$  Set $A_{\pi(n+1)}=\emptyset$,
$h_i=\dot{g}_{\pi(1)}^{-1}.\dot{g}_{\pi(2)}^{-1}.\dots.\dot{g}_{\pi(i)}^{-1}$ and $R_i=h_i(\tilde{B}_{\pi(i)}\setminus \tilde{A}_{\pi(i+1)})$,
for  $1\leq i\leq p$. By induction on $p$ we see
$$\tilde{A}_{\pi(1)}=\bigsqcup_{i=1}^pR_i.$$
Since by assumption $\sum_{i=1}^p(W_i-V_i)=(1,1,\dots,1)$,  the equations\begin{center}
 $\bigsqcup_{C\in (B_{\pi(i-1)}\setminus A_{\pi(i)})}x_0(C)=h_i^{-1}R_i,\ \ \ \ \ (1\leq i\leq p)$\end{center}
define a paradoxical decomposition. On the other hand we have
\begin{align*}
\tau(G)\leq& |B_{\pi(p)}|+(\sum_{i=2}^n |B_{\pi(i-1)}\setminus
A_{\pi(i)}|)+1\\ =&(\sum_{i=1}^n
|B_{\pi(i)}-A_{\pi(i)}|)+|A_{\pi(1)}|+1\\ =& \ell+|A_{\pi(1)}|+1.
\end{align*}
\end{proof}
In the special case where $G$ does not contain the free group on
two generators and $\ell=3$ and $|A_\pi(1)|=1$ we have $\tau(G)=5.$

\begin{example}\em \label{t5}
Let $\c=\{(1,2,2,2), (2,1,2,1), (2,2,1,1)\}$. Then $\tau(G)\leq 5$.\\
A minimal diagram associated with $Eq(\frak g,\mathcal E)$ is
 {\tiny
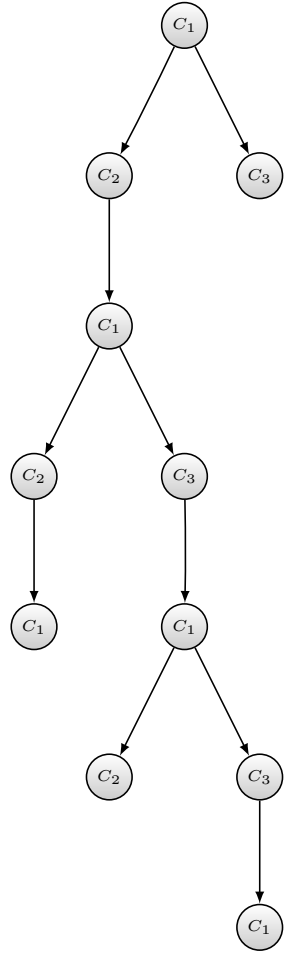
{
\begin{center}
\begin{tikzpicture}[-latex ,auto ,node distance =2 cm and 1cm ,on grid ,
semithick , state/.style ={ circle ,top color =white , bottom
color = black!20 , draw,black , text=black, minimum width =.3cm}]

\node[state] (C)[black] {$C_3$};
\node[state] (A) [black,above left=of C] {$C_1$};
\node[state] (B) [below left =of A] {$C_2$};
\node[state] (M) [below =of B] {$C_1$};
\node[state] (N) [below left=of M] {$C_2$};
\node[state] (G) [below right=of M] {$C_3$};
\node[state] (E) [below  =of G] {$C_1$};
\node[state] (E') [black, below  =of N] {$C_1$};
\node[state] (N') [black,below left=of E] {$C_2$};
\node[state] (G') [below right=of E] {$C_3$};
\node[state] (O) [black,below  =of G'] {$C_1$};

\path (A) edge [bend left =1] node[below =0.15 cm] {$$} (B);
\path (A) edge[bend right =0] node[below =0.15 cm] {$$} (C);
\path (M) edge[bend right =0] node[below =0.15 cm] {$$} (N);
\path (M) edge[bend right =0] node[below =0.15 cm] {$$} (G);
\path (G) edge [bend left =1] node[below =0.15 cm] {$$} (E);
\path (E) edge[bend right =0] node[below =0.15 cm] {$$} (N');
\path (E) edge[bend right =0] node[below =0.15 cm] {$$} (G');
\path (G') edge[bend right =0] node[below =0.15 cm] {$$} (O);
\path (N) edge[bend right =0] node[below =0.15 cm] {$$} (E');
\path (B) edge[bend right =0] node[below =0.15 cm] {$$} (M);

\end{tikzpicture}
\end{center}
}}
The number of paths from the top to the bottom is 4. These paths are $[C_1,C_3], $$[C_1,C_2,C_1,C_2,C_1],$ $[C_1,C_2,C_1,C_3,C_1,C_2]$ and $[C_1,C_2,C_1,C_3,C_1,C_3,C_1].$ Therefore $\tau(G)\leq 4+1=5$.   Indeed according to (\ref{key}) we have
\begin{align*}
  x_0(C_1) =& g_3x_0(C_2)\cup g_3x_0(C_3) \\
  x_0(C_2) =& g_1^{-1}x_0(C_1) \\
  x_0(C_3) =& g_2x_0(C_1).
\end{align*}
Define
\begin{align*}
  \Delta_1 =& g_3g_1^{-1}g_3g_1^{-1}g_3x_0(C_2), \\
  \Delta_2 =& g_3g_1^{-1}g_3g_1^{-1}g_3g_2x_0(C_1,) \\
  \Delta_3 =& g_3g_1^{-1}g_3g_2x_0(C_1), \\
  \Delta_4 =& g_3x_0(C_3), \\
  \Delta_5 =& x_0(C_2)\cup x_0(C_3).
\end{align*}
Then $x_0(C_1)=\bigsqcup_{i=1}^{4} \Delta_i$. So, $G=\bigsqcup_{i=1}^{5} \Delta_i$ and finally we have
\begin{align*}
  G= & g_2^{-1}g_3^{-1}g_1^{}g_3^{-1}\Delta_3\bigsqcup g_3^{-1}g_1g_3^{-1}g_1g_3^{-1}\Delta_1\bigsqcup g_3^{-1}\Delta_4\\
  G= & g_2^{-1}g_3^{-1}g_1g_3^{-1}g_1g_3^{-1}\Delta_2\bigsqcup e\Delta_5,
\end{align*}
where $e$ is the identity element of $G$. We have a complete paradoxical decomposition with five pieces. So, $\tau(G)\leq 5$.

\end{example}
\begin{example}\em \cite{r-w}
Let $\mathbb{F}_2=\langle a, b\rangle$ be the free group with generators $a$ and $b$ and identity element $e$. Suppose that $\frak g=(a,b)$ and $\mathcal E=\{E_1, E_2, E_3\}$,
where \\
\\
$E_1=\{x,\ \text{x is a reduced word starting with a}\}$,\\
$E_2=\{x,\ \text{x is a reduced word starting with b}\}$ and \\
$E_3=\{x,\ \text{x is a reduced word starting with}\ a^{-1}\ \text{or}\ b^{-1}\}\bigsqcup\{e\}.$\\

\noindent Then $Con(\frak g, \mathcal E)=\{C_1,\dots, C_7\}$, where $C_1=(1,1,2), C_2=(2,1,2), C_3=(3,1,2), C_4=(3,2,2), C_5=(3,3,2), C_6=(3,1,3)$  and
$C_7=(3,1,1).$ One can see that this example satisfies that condition of the pervious theorem.  Consider the following configuration equations
$$\sum_{x_0(C)\subseteq E_1}f_C=\sum_{x_1(C)\subseteq E_1}f_C\ \text{and}  \sum_{x_0(C)\subseteq E_2}f_C=\sum_{x_2(C)\subseteq E_2}f_C$$
i.e.
$$f_{C_1}=f_{C_1}+f_{C_2}+f_{C_3}+f_{C_6}+f_{C_7}\ \text{and}\ f_{C_2}=f_{C_1}+f_{C_2}+f_{C_3}+f_{C_4}+f_{C_5},$$
which imply that each $f_C$ is zero.
\noindent The corresponding diagram is

\vspace{.4cm}
{\tiny
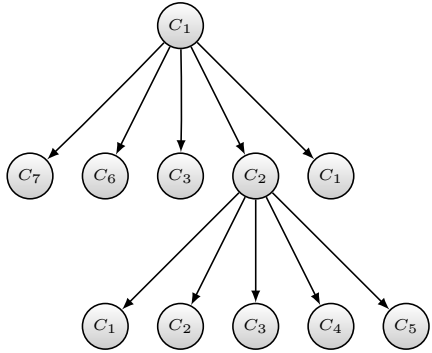
{
\begin{center}
\begin{tikzpicture}[-latex ,auto ,node distance =2 cm and 1cm ,on grid ,
semithick , state/.style ={ circle ,top color =white , bottom
color = black!20 , draw,black , text=black, minimum width =.3cm}]

\node[state] (C)[black] {$C_6$};
\node[state] (M) [right =of C] {$C_3$};
\node[state] (A) [black,above =of M] {$C_1$};
\node[state] (MM) [right =of M] {$C_2$};
\node[state] (MMM) [right =of MM] {$C_1$};
\node[state] (MMMM) [left =of C] {$C_7$};

\node[state] (AA) [black,below =of M] {$C_2$};
\node[state] (MA) [right =of AA] {$C_3$};
\node[state] (MMA) [right =of MA] {$C_4$};
\node[state] (MMMA) [right =of MMA] {$C_5$};
\node[state] (MMMMA) [left =of AA] {$C_1$};

\path (A) edge [bend left =1] node[below =0.15 cm] {$$} (M);
\path (A) edge[bend right =0] node[below =0.15 cm] {$$} (C);
\path (A) edge[bend right =0] node[below =0.15 cm] {$$} (MM);
\path (A) edge[bend right =0] node[below =0.15 cm] {$$} (MMM);
\path (A) edge [bend left =1] node[below =0.15 cm] {$$} (MMMM);
\path (MM) edge[bend right =0] node[below =0.15 cm] {$$} (AA);
\path (MM) edge[bend right =0] node[below =0.15 cm] {$$} (MA);
\path (MM) edge[bend right =0] node[below =0.15 cm] {$$} (MMA);
\path (MM) edge[bend right =0] node[below =0.15 cm] {$$} (MMMA);
\path (MM) edge[bend right =0] node[below =0.15 cm] {$$} (MMMMA);

\end{tikzpicture}
\end{center}
}}
\vspace{.8cm}
 \noindent Put $A=a^{-1}x_0(C_1),$
$B=b^{-1}a^{-1}[\cup_{i=1}^5 x_0(C_i)], \ C=a^{-1}(x_0(C_6)\cup x_0(C_7))$ and $M=\mathbb{F}_2\setminus x_0(C_1).$
Then we have the non-complete paradoxical decomposition $$\mathbb{F}_2=A\sqcup B\sqcup C \sqcup M \sqcup a^{-1} x_0(C_3)=eM\sqcup aA=abB\sqcup aC.$$
Hence $\tau(\mathbb{F}_2)= 4$.
\end{example}

We now summarize what we have done for estimating the upper bound of the Tarski number of a group.
\begin{enumerate}
  \item Change the order of the configuration equations based on the permutation matrix in the definition of the normality condition.
  \item Define the initial diagrams consisting of $p$ connected components.
  \item Make each set in the top of the above diagrams disjoint from the previous ones by borrowing their intersection from the sets in the bottom of the former diagrams.
  \item Glue these diagrams as it is explained in subsection 4.1.
  \item Find the configurations of the top sets of the initial diagrams in the bottom of the last gluing diagram.
  \item Count all the directed pats from top to the bottom. This number plus one is an upper bound for $\tau(G)$.
   \end{enumerate}
To get the best bound for $\tau(G)$ using this procedure, we borrow the maximum number of configurations that we can get from the former bottom sets (not only the intersections). It is applied in the examples.

\begin{definition}
Let $m,n\in \mathbb{N}.$ A set of $(n+1)$-tuples $\mathcal C=\{(c^i_0,\dots, c^i_n),\ \ 1\leq i\leq \ell\}$ with $1\leq c^i_j\leq m$ and  $$\cup_{j=1}^n \{c^i_j\}=\{1,\dots, m\},\ \ \  \  (i=1,\dots,\ell)$$  is called a pre-configuration set if there exist a group $G$, a string $\frak g$ of elements of $G$ and a partition $\mathcal E$ of $G$ such that $\c=\mathcal C$.
\end{definition}

In \cite{sapir ershov} the authors give examples of groups with Tasrki numbers 5 and 6.
Now the question is whether we can construct such groups using configurations.
Responding to this question depends on knowing that given well-behaved sets are pre-configuration ones.
In particular
\begin{question}
 Is $\mathcal C=\{(1,2,3,2), (1,3,1,3), (2,1,2,2), (3,3,1,2), (3,3,2,1)\}$ a pre-configuration set (see Example \ref{tt5})?
\end{question}


\end{document}